\documentclass{amsart}
\usepackage[utf8]{inputenc}
\usepackage{hyperref}
\usepackage{verbatim}
\usepackage{graphicx}
\usepackage{comment}
\usepackage{amsmath}
\usepackage{amsthm}
\usepackage{amssymb}
\usepackage{cleveref}
\usepackage{mathrsfs}
\usepackage{mathtools}
\usepackage{bm}
\newcommand{\real}{\mathbb{R}}
\newcommand{\complex}{\mathbb{C}}

\newcommand{\M}[1]{\left({#1}\right)}
\newcommand{\Mb}[1]{\left[{#1}\right]}

\newcommand{\rmi}{{\mathrm{ i}}}

\newcommand{\sfu}{\mathsf{u}}

\renewcommand{\Re}{\textrm{Re}}
\renewcommand{\Im}{\textrm{Im}}

\newcommand{\captionfont}{}

\newtheorem{theorem}{\bf Theorem}

\newtheorem{lem}{\bf Lemma}
\newtheorem{Rem}{\bf Remark}
 
\begin{document}
\title[Active exterior cloaking for Helmholtz eq. with complex wavenumbers]{Active exterior cloaking for the 2D Helmholtz equation with complex wavenumbers and application to thermal cloaking}
\author[M. Cassier, T. DeGiovanni, S. Guenneau, F. Guevara Vasquez]{Maxence Cassier$^{1}$, Trent DeGiovanni$^{2}$, S\'ebastien Guenneau$^{3}$, and Fernando Guevara Vasquez$^{2}$}
\address{$^{1}$Aix Marseille Univ, CNRS, Centrale Marseille, Institut Fresnel, Marseille, France\\
$^{2}$University of Utah, Mathematics Department, Salt Lake City UT 84112, USA\\
$^{3}$UMI 2004 Abraham de Moivre-CNRS, Imperial College London, London SW7 2AZ, UK
}


\begin{abstract}
    We design sources for the two-dimensional Helmholtz equation that can cloak an object by cancelling out the incident field in a region, without the sources completely surrounding the object to hide. As in previous work for real positive wavenumbers, the sources are also determined by the Green identities. The novelty is that we prove that the same approach works for complex wavenumbers which makes it applicable to a variety of media, including media with dispersion, loss and gain. 
    Furthermore, by deriving bounds on Graf's addition formulas with complex arguments, we obtain new estimates that allow to quantify the quality of the cloaking effect.
    We illustrate our results by applying them to achieve active exterior cloaking for the heat equation. 
\end{abstract}
\maketitle

\section{Introduction}
Our goal is to use specially designed sources to cloak or hide a bounded object from a probing field $u_i$ satisfying the two dimensional Helmholtz equation 
\begin{equation}
    \Delta u_i + k^2 u_i = 0,
    \label{eq:helmholtz}
\end{equation} 
in a region containing the object. Here $\Delta$ denotes the Laplacian. This is called \textit{active cloaking} since to build the cloak, we use sources rather than passive materials that may be hard to manufacture \cite{kadic2013metamaterials}. Moreover, a great
\noindent advantage of the cloaking strategy we present is that it does not require completely 
surrounding the object to hide it, hence the \textit{exterior cloaking} name. This idea was introduced in \cite{Guevara:2009:AEC} for the 2D Helmholtz equation with $k$ real (lossless propagative media). Here we allow $k$ to take any values in the complex plane, except for the negative real axis. Thus using a frequency decomposition of the transient regime via a Fourier-Laplace transform on the time variable, our approach applies to cloaking objects for acoustic waves propagating in passive, dissipative, active, or dispersive media, and similarly for diffusive media. Interestingly, complex wavenumbers open a path to exterior cloaking for problems modelled by partial differential equations with second order derivatives in space and with time derivatives of an arbitrary order. Moreover, we derive new error estimates on the convergence of active exterior cloaking and apply our results to cloaking for the heat equation in the transient regime.

\subsection{Active exterior cloaking}
\label{sec:aec_intro}
From potential theory \cite{Kellogg:1967:FPT} or using the Green identities (see e.g. \cite{Colton:2013:IAE}), it is possible to reproduce a solution to the homogeneous Helmholtz equation inside of a bounded region $\Omega$ and getting, simultaneously, a zero field outside of $\Omega$. This is achieved by a distribution of monopole and dipole sources on the boundary $\partial \Omega$ that can be expressed in terms of the value of the field and its normal derivative on $\partial \Omega$. As observed by Miller \cite{Miller:2005:OPC}, this principle can be used for cloaking. Indeed the monopole and dipole distribution can be chosen to generate the \emph{cloak field}
\begin{equation}
 u_c = \begin{cases} -u_i & \text{in } \Omega\\
 0 & \text{outside }\overline{\Omega},
 \end{cases}
 \label{eq:ucfield}
\end{equation}
where $\overline{X}$ denotes the closure of a set $X$.
We see that by linearity, $u_c + u_i$ cancels out inside $\Omega$ without affecting $u_i$ outside of $\Omega$. The end result is that objects inside $\Omega$ will not scatter and it is impossible to detect the cloaked field outside of $\Omega$. We call this approach the \emph{Green identity cloak}. A first drawback of this approach
is that the probing field $u_i$ needs to be known ahead of time. A second drawback is that the sources completely surround the object that we wish to hide. The exterior cloaking approach lifts this second limitation.

\begin{figure}
     \centering
     \includegraphics[width=0.5\textwidth]{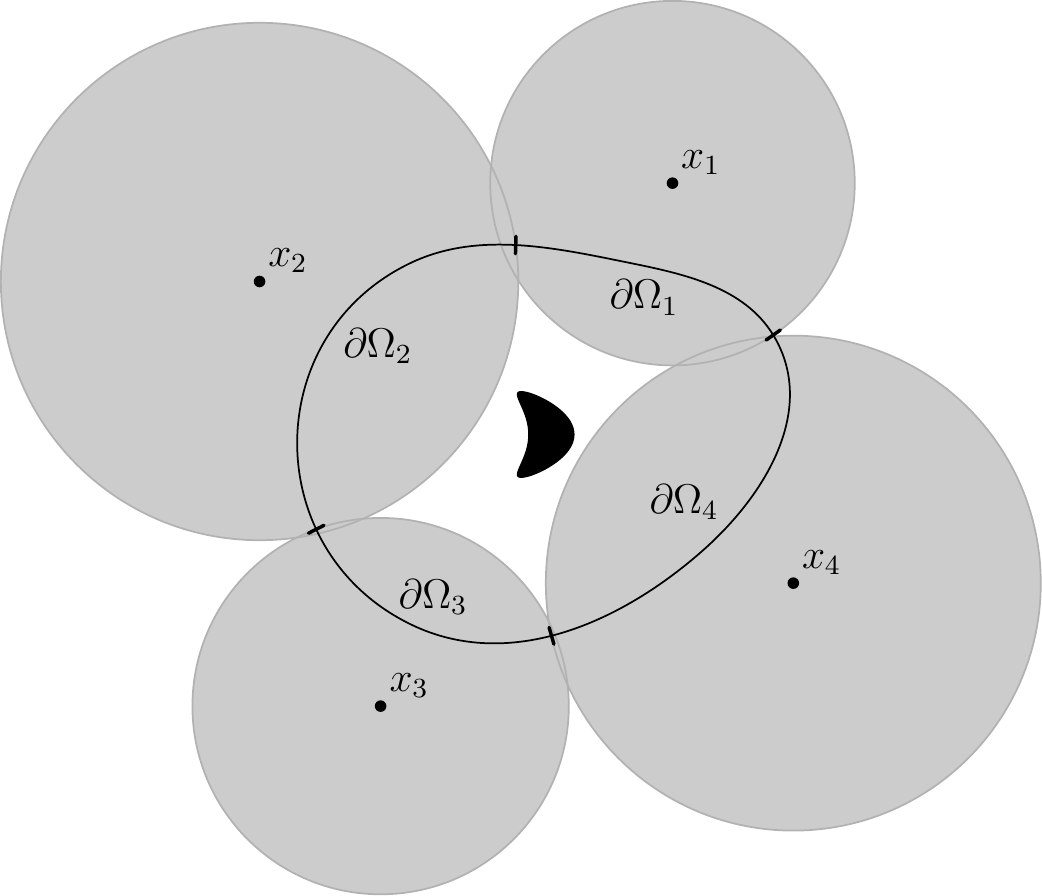} 
     
    \caption{Active exterior cloak for the Helmholtz equation starting from the Green identities applied on the surface $\partial\Omega$ of a domain $\Omega$. To hide the object (kite) inside $\Omega$, the sources in the portions $\partial\Omega_j$ of $\partial\Omega$ are moved to new locations $x_j$, $j=1,\ldots,N_{dev}$. The gray disks are the domains of divergence $R_j$ for the new fields, see \eqref{eq:rj}.  For illustration purposes we took $N_{dev}=4$, but 3 sources would suffice to achieve exterior cloaking.}
    \label{fig:setup}
\end{figure}

To achieve exterior cloaking we follow the approach in \cite{Guevara:2011:ECA} for the 2D Helmholtz equation with real wavenumbers, see also \cite{Guevara:2013:TEA} for the 3D Helmholtz equation and \cite{Norris:2014:AEC,Oneil:2015:ACI} for elasticity. See also \cite{Onofrei:2014:AMF} for a general analysis. The key observation is that Graf's addition formulas (see e.g. \S 10.23 in \cite{NIST:DLMF} or \cite{Martin:2006:MS}) can be used to move a monopole (or dipole) located at $y$ to a new location $x_j$. However the price to pay is that the new source is obtained by an infinite superposition of multipolar sources that diverges in the disk 
\begin{equation}\label{eq.domaindivergence}
D_{x_j,y}:=\{ x \in \mathbb{R}^2 \mid  |x-x_j| \leq |y-x_j|\},
\end{equation}
where $|\,\cdot\,|$ denotes the Euclidean norm. By linearity, we should also be able to move a distribution of monopoles or dipoles on a compact portion $\partial \Omega_j$ of the boundary $\partial \Omega$ to a new location $x_j$, obtaining the same cloak field $u_c$, provided we are outside of the closed disk
\begin{equation}
     \label{eq:rj}
     R_j = \{ x \in \mathbb{R}^2 \mid  |x-x_j| \leq \max_{y\in \partial\Omega_j} |y-x_j|\}.
\end{equation}
Assuming the portions $\partial \Omega_j$ cover $\partial \Omega$ and their intersection is reduced to points, we can achieve the same cloaking effect as the Green identity cloak if we are outside of the region $R_1 \cup R_2 \cup \ldots \cup R_{N_{dev}}$. Here we extend this approach to complex wavenumbers which enables new applications. Moreover, by obtaining bounds for the Graf's addition formula with complex arguments, we derive a simple geometric series ansatz to predict the truncation error for the field generated by the multipolar sources. This extends previous work on the truncation error of Graf's addition formulae \cite{Cassier:2013:MSA,Meng:2016:BTE} for real arguments. The truncation estimates allow us to quantitatively predict the quality of the cloaking effect.

\subsection{Extending active exterior cloaking to complex $k$}
Many partial differential equations in the frequency domain lead to the Helmholtz equation with a complex wavenumber. To name a few: the telegraph equation, the diffusion equation, Schr\"odinger equation and the Klein-Gordon equation (see e.g. \cite{Dautray:1992:MAN,ockendon2003applied} or \cite[\S 1.1.2]{gumerov2005fast}).
More generally, consider partial differential equations (PDEs) in the time domain of the form
\begin{equation}
 P(\partial_t) \mathfrak{u} = \Delta \mathfrak{u} + \mathfrak{f},
 \label{eq:poly}
\end{equation}
where $P$ is a polynomial of degree $n$ and $\mathfrak{f}(x,t)$ is a source term. Since \eqref{eq:poly} is a constant coefficient PDE, it admits a solution in the distributional sense, for example for any compactly supported source term $\mathfrak{f}$. This can be seen from the Malgrange-Ehrenpreis theorem, though without uniqueness or causality guarantees, see e.g. \cite{Friedlander:1998:ITD}.

Many classic equations are of the form \eqref{eq:poly}. For example, the wave equation can be obtained with $P(z) = z^2$, and the heat equation with $P(z) = z$. For a \emph{causal} source $\mathfrak{f}$ (i.e. $\mathfrak{f}(x,t) = 0$ for $t<0$), we can analyze equations of the form \eqref{eq:poly} in the frequency domain by means of the Fourier-Laplace transform
\begin{equation}\label{eqn:lp_transform}
u(x,\omega)= \int_0^\infty dt [e^{\rmi\omega t} \mathfrak{u}(x,t)],
\end{equation} 
where $\omega$ is in general complex and $\mathfrak{u}$ is assumed to grow sufficiently slowly. For example, we may assume that $\mathfrak{u}(x,t)$ satisfies for $t>0$
\begin{equation}
     \| \mathfrak{u}(\cdot,t)\| \leq C e^{\alpha t} (1+t^p),
     \label{eq:growth}
\end{equation} 
where $C>0$, $p \in \mathbb{N}$, and $\alpha \in \real$ are constants,  see e.g. \cite{arendt2001cauchy,Dautray:1992:MAN}. Under this assumption, the Fourier-Laplace transform is defined on the half plane
\begin{equation}
\mathbb{C}_{\alpha}^+:=\{ \omega \in \complex ~|~ \Im(\omega) > \alpha\}.
\label{eq:calpha}
\end{equation}
The  choice of norm depends on the spatial differential operator. Since we focus on the Laplacian, we use the $H^1$ norm on any bounded open set of interest ($\alpha$, $p$, $C$ could depend on the choice of the set). To summarize, in our situation we may assume that $\mathfrak{u}$ satisfies the growth condition \eqref{eq:growth} and $\mathfrak{u}\in L^1_{loc}( (0,\infty), H^1_{loc}(\real^2) )$ \footnote{We recall that for $p\geq 1$, $L^p_{loc}(X)$ (resp.  $H^1_{loc}$) is the set of functions that are $L^p$ (resp. $H^1$) on any bounded open subset of $X$, with closure inside $X$. See e.g. \cite{Evans:2010:PDE,Lions:1972:NHB}.} and all its time derivatives up to order $n$ (the degree of the polynomial $P$) and the source term $\mathfrak{f}$ are in $L^1_{loc}( (0,\infty), L^2_{loc}(\real^2) )$ and satisfy also \eqref{eq:growth}. This growth condition allows to make sense of the Fourier-Laplace transform for solutions that may grow exponentially in time\footnote{The growth condition \eqref{eqn:lp_transform} and the regularity assumption $\mathfrak{u}\in L^1_{loc}( (0,\infty), L^2_{loc}(\real^2) )$  ensure in particular the existence of the Laplace transform  \eqref{eqn:lp_transform} as a Bochner integral with respect to $t$ valued in $L^2$ \cite{arendt2001cauchy}, and thus also pointwise for a. e. $x\in \mathbb{R}^2$ and all $\omega\in \mathbb{C}^+_{\alpha}$. We point out that in a more general context, the Fourier-Laplace transform can be extended to spaces of distributions, see e.g. \cite{Dautray:1992:MAN,Sayas:2016:RPT,Zemanian:1972:RTC}.}, as is the case of active media\footnote{By ``active media'' we mean there is energy input that may lead to increase of the magnitude of the fields in time, see for instance \cite{Skaar:2006:Fresnel}. By ``gain media'' we mean that there are spatially growing outgoing solutions of the Helmholtz equation as for e.g. resonant states, see \cite{Dyatlov:2019:MTS}. This is not a universal nomenclature.}. 

\begin{figure}
 \centering
 \includegraphics[width=0.4\textwidth]{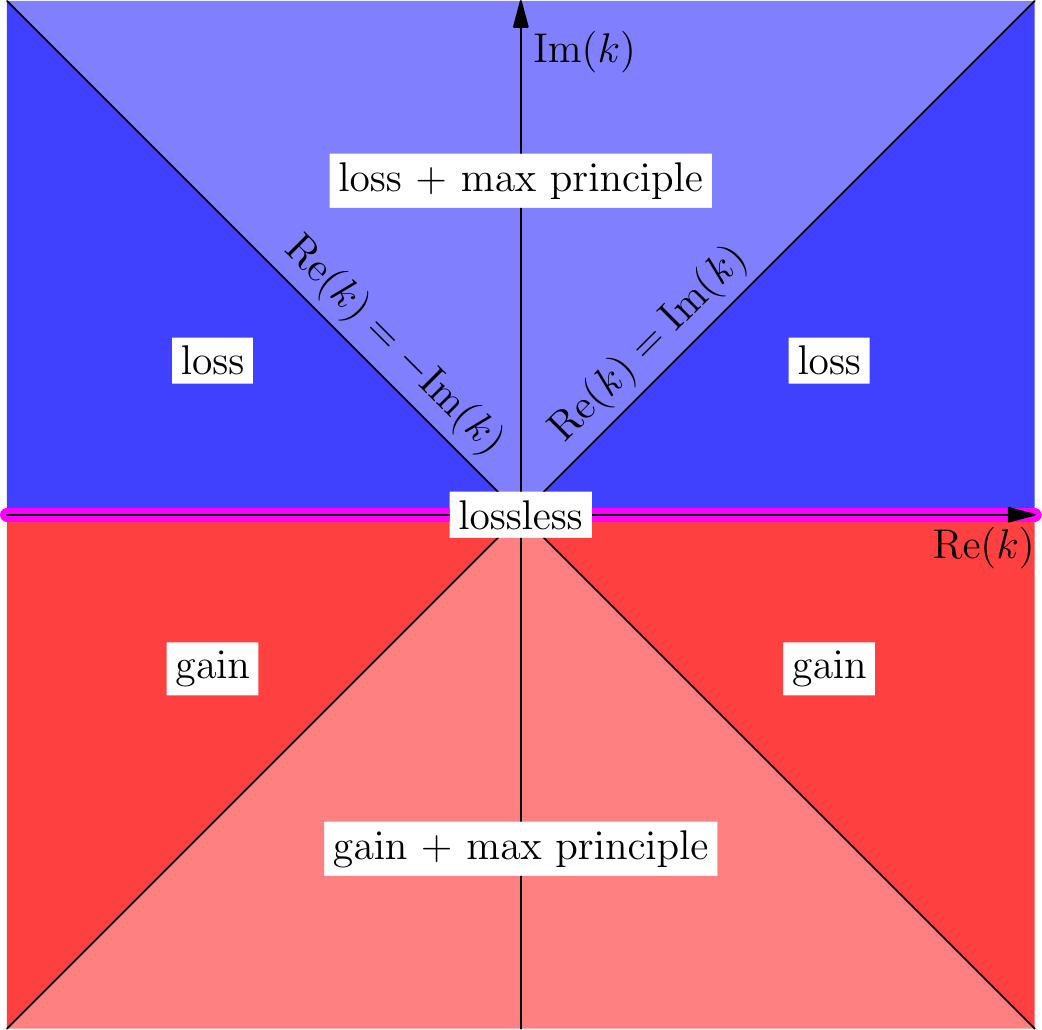}
 \caption{Diagram of cases for the wavenumber $k$ defined by  $k = \pm \rmi \sqrt{P(-\rmi\omega)}$, where we used the principal square root and the sign is chosen to match the sign of $\Im(k)$.
 }
 \label{fig:kdiagram}
\end{figure}

Assuming $\mathfrak{u}$ and all its time derivatives up to order $n-1$ vanish at $t=0$, we get that $u(x,\omega)$ satisfies the Helmholtz equation:
\begin{equation}
    \Delta u + k^2 u = -f,
    \label{eq:helmholtz2}
\end{equation}
with the relation $ 
 k^2 = -P(-\rmi\omega).
$
This formalism shows that $k$ may be complex e.g. in the case of the heat equation. Another example is the modified wave equation with $P(z) = z^2 + \alpha z$, $\alpha \in \real$. If $\alpha>0$, this is the dissipative wave equation which corresponds to wave propagation in lossy media. If $\omega>0$, we can choose the root $k$ such that $\Im(k)>0$ to get spatially decreasing solutions to \eqref{eq:helmholtz2}\footnote{We shall see in \cref{sec:conv} that the decay is consistent with the choice of Green function \eqref{eq:greenfun}.}. Whereas with $\alpha<0$, $\omega>0$ we choose $k$ such that $\Im(k)<0$ to get spatially increasing solutions corresponding to an amplifying medium (medium with gain), see e.g. \cite{Dyatlov:2019:MTS}.
We summarize the possibilities in \cref{fig:kdiagram}. We have highlighted  the region with $\Re (k^2) < 0$ (or equivalently $|\Re(k)| < |\Im(k)|$), since a form of the maximum principle holds for the Helmholtz equation with such $k$ \cite{Kresin:1993,Krutitskii:1998:TDP}. Later in \cref{sec:stability}, we see that the maximum principle gives a form of stability for the accuracy of approximations to our cloaking approach.

Other situations where complex wavenumber $k$ arises are in passive, dispersive media, where the index of refraction is a complex valued function of frequency $\omega$ \cite{Tip:1998,Figotin:2005,Gralak:2010,Cassier:2017,Cassier:2017bis,Bellis:2019}. A typical example is the dispersion law given by the Drude-Lorentz model in electromagnetism, which can model both metals and metamaterials with negative index of refraction \cite{ramakrishna2005physics}.  In acoustics (see for e.g. section 1.1.2 of \cite{gumerov2005fast}), complex wavenumbers also arise when studying acoustic (pressure) waves propagating within complex (but homogeneous and isotropic) fluids (i.e. not barotropic) with a relaxation time (due to the presence of solid particles or bubbles) which can be modelled in the time-harmonic regime with $k^2=\omega^2/(c^2(1-\rmi\omega\tau_p))$. Here $c$ is the speed of sound (m/s), $\omega$ the pulsation frequency (rad/s) and $\tau_p$ the density relaxation time (s). 

We note that the frequency domain formulation gives a strategy for active exterior cloaking in the time domain for equations of the form \eqref{eq:poly}, and even when the powers are fractional or negative (which corresponds to integro-differential equations). One caveat of our approach is that the sources need to make sense physically in \eqref{eq:poly}. However, for active thermal cloaking (\cref{sec:timedomain}) the sources we obtain can be thought of as Peltier devices \cite{DiSalvo:1999:TCP}.

\subsection{Structure of the paper}

We derive convergence estimates for the Graf addition formula applied to Green functions in \cref{sec:conv}, showing that the truncation error of the series can be dominated by that of a geometric series with ratio that depends only on the position of the evaluation point relative to the positions of the original and new sources. This result can be applied to get truncation estimates for the multipolar source expansions that appear in active exterior cloaking for possibly complex wavenumbers (\cref{sec:aec}). In \cref{sec:aec}, we also use a form of the maximum principle for the Helmholtz equation, which guarantees the truncation errors in a region are maximum on the boundary of the region (this only holds for  a class of dissipative or diffusive media). The time domain problem for the heat equation is then considered in \cref{sec:timedomain}. We conclude with some future work and perspectives in \cref{sec:discussion}.

\section{Moving sources}\label{sec:conv}
The field evaluated at $x$ corresponding to a point source located at $y$ is given by the appropriate Green function
\begin{equation}
 G(x-y;k) = \frac{\rmi}{4} H_0^{(1)}(k|x-y|),
 \label{eq:greenfun}
\end{equation}
where $H_0^{(1)}$ is the zero-th order Hankel function of the first kind\footnote{In the context of the wave equation with constant propagation speed $c$, we have $k = \omega/c >0$ and the choice of Green function \eqref{eq:greenfun} corresponds to outgoing waves. This is consistent with the Fourier-Laplace transform convention \eqref{eqn:lp_transform} and with the convention that the corresponding time harmonic field is $\Re(\exp[-\rmi\omega t] G(x-y;k))$. In fact $\Re(k)>0$ also gives outgoing waves, as can be seen e.g. from adapting the discussion \cite[eq. 1.2.12]{gumerov2005fast} from 3D to 2D using  \cite[eq. 10.4.3]{NIST:DLMF}.}, see e.g. \cite[eq. 10.4.3]{NIST:DLMF}. 
Moreover $G(x;k) \to 0$ as $|x| \to \infty$ whenever $\Im(k) \geq 0$, as can be seen from the large argument asymptotic for Bessel functions \cite[eq. 10.2.5]{NIST:DLMF}. 
When $\Im(k) <0$, the same asymptotic shows that $|G(x;k)| \to \infty$. Thus the choice of Green function is consistent with the loss and gain conventions in the diagram appearing in \cref{fig:kdiagram}.

Thanks to the Graf addition formulas \cite[eq. 10.23.7]{NIST:DLMF}, we can move (with three significant caveats) the source from location $y$ to another location $x_j$, indeed:
\begin{equation}\label{eqn:grafs}
G(x-y;k) = \frac{\rmi}{4}\sum_{m=-\infty}^\infty H_m^{(1)}(k|x-x_j|)J_m(k|y-x_j|)\exp[\rmi m\theta], 
\end{equation}
where $\theta = \arg(x-x_j)-\arg(y-x_j)$ and $\arg x$ is the counter-clockwise angle between the vectors $x$ and $(1,0)$. Here $J_m$ is the $m-$th order Bessel function of the first kind, see e.g. \cite[eq. 10.2.2]{NIST:DLMF}. The first caveat is that the new source in \eqref{eqn:grafs} is no longer a monopole point source like \eqref{eq:greenfun}, but a linear combination of Helmholtz equation solutions that diverge as $|x-x_j| \to 0$, of the form $V_m(x-x_j)$ where 
\begin{equation} 
 V_m (x) = \exp[\rmi m\arg(x)] H_m^{(1)}(k|x|),
 \label{eq:vm}
\end{equation}
and that are known as \emph{multipolar sources} (or \emph{cylindrical outgoing waves} when $k$ is real).
The second caveat is that the Graf addition formula is only valid for $k \in \complex \setminus (-\infty,0]$. The third caveat is that the series converges only outside of the disk $D_{x_j,y}$, as defined in \eqref{eq.domaindivergence}. 

The same method and caveats apply if we desire to move a dipole located at $y$ and oriented in the direction $\nu(y)$ normal to the boundary $\partial \Omega$ at $y$, or more precisely
\begin{equation}\label{eqn:grafs.dipole}
\frac{\partial G}{\partial \nu(y)}(x-y;k)=\frac{\rmi}{4}\sum_{m=-\infty}^\infty H_m^{(1)}(k|x-x_j|) \frac{\partial }{\partial \nu(y)}\big(J_m(k|y-x_j|) \exp[\rmi m\theta]\big),
\end{equation}
where $\theta$ is the same as in \eqref{eqn:grafs}. Formally speaking, equation \eqref{eqn:grafs.dipole} can be obtained by taking the gradient term by term (with respect to $y$) in \eqref{eqn:grafs}  and then taking the dot product with $\nu(y)$. 
The differentiation term by term can be easily justified by using \cref{lem:gbound} in order to prove that the involved series of gradients is locally normally convergent (and thus locally uniformly convergent) with respect to $y$ when $x$, $x_j$ are fixed.

We start in \cref{sec:tee} by proving convergence estimates for \eqref{eqn:grafs} and \eqref{eqn:grafs.dipole}. The convergence errors are illustrated numerically in \cref{sec:ten}.

\subsection{Truncation error estimates}
\label{sec:tee}

To study the convergence rate of \eqref{eqn:grafs} and \eqref{eqn:grafs.dipole} we define the truncation to $2M+1$ terms of the formula \eqref{eqn:grafs} for moving the point source at $x_j$ to location $y$ by
\begin{equation}\label{eqn:grafstr}
G_{j,M}(x-y;k) = \frac{\rmi}{4}\sum_{m=-M}^M H_m^{(1)}(k|x-x_j|)J_m(k|y-x_j|)\exp[\rmi m\theta], 
\end{equation}
where $M$ is an integer. The truncation error for a monopole and dipole are given respectively by  
\begin{equation}
\begin{aligned}
R_{j,M}(x;k) &= |G(x-y;k) - G_{j,M}(x-y;k)| \mbox{ and }\\
R'_{j,M}(x;k)&= \left|\frac{\partial}{\partial \nu(y)} [G(x-y;k)- G_{j,M}(x-y;k) ]\right|.
\end{aligned}
\label{eq:rjm}
\end{equation}
In the next theorem, we show that these truncation errors are dominated by the truncation errors of  well-known series such as geometric series. Our convergence estimates account for moving sources from different original positions $y$ to a single new position $x_j$. The case of different original positions $y$ is useful in the context of active exterior cloaking (\cref{sec:aec}).
We point out that the monopole truncation error was derived for real wavenumbers by \cite[Lemma 9]{Cassier:2013:MSA}  and \cite{Meng:2016:BTE}, using techniques that are similar to the ones we use here. \Cref{pro:grafs_trunc} applies to the  monopole and the dipole truncation errors and allows for estimates that are \emph{uniform} with respect to original source location $y$, evaluation point $x$ and complex wavenumbers $k$.
 
\begin{theorem}
\label{pro:grafs_trunc}
Let $x_j \in \real^2$, $M\geq 2$, $Y\subset \real^2$ be a compact set and define the disk 
\[D_j^{\max} =  \{ x \in \mathbb{R}^2 ~\big\vert~ |x-x_j| \leq \max_{y \in Y}|y-x_j| \}.\]
Let $X \subset \real^2 \setminus D_j^{\max}$ and $K \subset \complex \setminus(-\infty,0]$ be compact sets. Then for any $(x,k) \in X \times K$ we have the following bounds for the monopole and dipole truncation errors
\begin{equation}\label{eq.graf}
\begin{split}
        R_{j,M}(x;k) & \leq C_1 \Big(-\ln(1-a_x)-\sum_{m=1}^M \frac{a_x^m}{m} \Big)
        \\ 
        R'_{j,M}(x;k) &\leq  C_2\frac{a_x^{M+1}}{1-a_x},
\end{split}
\end{equation}
where $C_1$ and $C_2$ may depend on $X$, $Y$, $K$ and
\begin{equation}
a_x :=  \frac{\max_{y\in Y} |y-x_j|}{|x-x_j|}.
\end{equation}
\end{theorem}

To prove \cref{pro:grafs_trunc} we need the following asymptotic formulas for Bessel functions that are uniform on the order and are valid on appropriate compact sets of $\complex$ excluding the negative real axis $(-\infty,0]$. This is because we use the power series definitions for Bessel functions in \cite[\S 10.8]{NIST:DLMF} and the power series expansion for $H_n^{(1)}(z)$ is not valid for $z \in (-\infty,0]$ (as the expansion contains the term $(2/\pi) \ln(z/2) J_n(z)$ which has a discontinuity for such $z$).
\begin{lem}\label{lem:gbound}
Let $K_1$ be a compact set of $\mathbb{C}$, $K_2$ be a compact set of $\mathbb{C}\setminus (-\infty,0]$ and $n\in \mathbb{N}$. Then there exists constant $C_{K_1}$, $\tilde{C}_{K_1}$ and $\tilde{C}_{K_2}$ (independent of $n$) such that:
\begin{eqnarray}
&\Big|J_n(z)-\displaystyle \frac{1}{n!}\Big(\frac{z}{2}\Big)^{n}\Big|&\leq \frac{C_{K_1}} {(n+1)!}\Big( \frac{|z|}{2}\Big)^{n+2},  \quad \forall z\in K_1, \ n\geq 0, \label{eq.Jn}\\
&\Big|J_n'(z) -\displaystyle \frac{1}{2 (n-1)!}\big(\frac{z}{2}\big)^{n-1}\displaystyle  \Big|&\leq   \frac{ \tilde{C}_{K_1}} {n!}\Big( \frac{|z|}{2}\Big)^{n+1}, \  \forall z\in K_1, \ n\geq 1, \label{eq.Jnprime}\\
 &\Big|H_n^{(1)}(z)-\displaystyle  \frac{\rmi  \, (n-1)! }{\pi}\Big(\frac{2}{z}\Big)^n  \displaystyle \Big|&\leq \tilde{C}_{K_2} (n-2)!  \Big(\frac{2 }{|z|}\Big)^{n-2},  \, \forall z\in K_2, \, n \geq 2.   \ \ \label{eq.Hn}
 \end{eqnarray}
\end{lem}

The proof of \cref{lem:gbound} is included in \cref{app:gbound}. From \cref{lem:gbound}, we can deduce the following inequalities for $J_n$, $J'_n$  and $H_n^{(1)}$ that are useful in the proof of \cref{pro:grafs_trunc}. 
Let $K_1$ be a compact set of $\mathbb{C}$ and $K_2$ be a compact set of $\mathbb{C}\setminus (-\infty,0] $. Applying the inequality \eqref{eq.Jn}, we get that for all $z\in K_1$ and $n\geq 0$ 
\begin{eqnarray}\label{eq.Jnbis-rem}
|J_n(z)|&\leq &\big|J_n(z)-\displaystyle \frac{1}{n!}\Big(\frac{z}{2}\Big)^{n}|+\displaystyle \frac{1}{n!}\Big(\frac{|z|}{2}\Big)^{n} \nonumber\\ 
&\leq& \frac{C_{K_1}} {(n+1)!}\Big( \frac{|z|}{2}\Big)^{n+2}+\displaystyle \frac{1}{n!}\Big(\frac{|z|}{2}\Big)^{n}  \nonumber \\
&\leq& \frac{B_{K_1}}{n!}\Big(\frac{|z|}{2}\Big)^{n} \ \mbox{ with } \ B_{K_1}=\max\Big(1, C_{K_1}\max_{z\in K_1}\Big(\frac{|z|}{2}\Big)^2\Big)>0.
\end{eqnarray}
Similarly, one deduces from formula \eqref{eq.Jnprime} and \eqref{eq.Hn} that there exists two constants  $\tilde{B}_{K_1}>0$ and $\tilde{B}_{K_2}>0$ such that:
\begin{eqnarray}
&\Big|J_n'(z) \displaystyle  \Big|\leq \displaystyle \frac{\tilde{B}_{K_1}}{ (n-1)!}\Big(\frac{|z|}{2}\Big)^{n-1}, &\,  \forall z\in K_1, \ n\geq 1, \label{eq.Jnprimebis-rem}\\
&\Big|H_n^{(1)}(z) \Big|\leq  \displaystyle \tilde{B}_{K_2} (n-1)!\Big(\frac{2}{|z|}\Big)^n, & \, \forall z\in K_2, \, n \geq 2. \label{eq.Hn-rem}
\end{eqnarray}
We remark that \cref{lem:gbound} shows that the inequalities \eqref{eq.Jnbis-rem}, \eqref{eq.Jnprimebis-rem} and \eqref{eq.Hn-rem} are optimal in the sense that they bound the functions by their leading order term.
We are now ready to prove \cref{pro:grafs_trunc}.

\begin{proof}
{\bf Step 1: inequality on the monopole truncation error $R_{j,M}$.\\}
 We want to apply the \cref{lem:gbound} to bound the terms $J_m(k|y-x_j|)$ and  $H^{(1)}_m(k|x-x_j|)$  for $m\geq M+1$ appearing in the expression of  $R_{j,M} = |G(x-y;k) - G_M(x-y;k)|.$ 
 Noting first that $H^{(1)}_{-m}=(-1)^mH^{(1)}_m$ and $J_{-m}=(-1)^mJ_m$,
 we obtain that:
\begin{align*}
R_{j,M}(x;k) \leq \sum_{|m|\geq M+1}|H_m^{(1)}(k|x-x_j|)J_m(k|y-x_j|)|
&=\sum_{m=M+1}^\infty 2|H_m^{(1)}(k|x-x_j|)J_m(k|y-x_j|)|.
\end{align*}
Thus, applying inequalities \eqref{eq.Jnbis-rem} and \eqref{eq.Jnprimebis-rem}  gives that there exists $C_1>0$ (depending on the compacts $X$ and $K$ but not on the truncation index $M$) such that
$$
R_{j,M}(x;k)\leq C_1 \displaystyle \sum_{m=M+1}^{\infty} \frac{1}{m} a_x^m=C_1 \bigg(-\ln(1-a_x)-\sum_{m=1}^M \frac{a_x^m}{m} \bigg)
.$$

\noindent {\bf Step 2: inequality on the dipole truncation error $R'_{j,M}$.\\}
For  $R'_{j,M}(x;k)$, we have the following (note that $\theta$ and $n$ depend on $y$) 
\begin{align*}
    R'_{j,M}(x;k)&= \bigg|\sum_{|m|\geq M+1} H_m^{(1)}(k|x-x_j|)\frac{\partial}{\partial \nu} J_m(k|y-x_j|)\exp[\rmi m\theta] \bigg| \\
     &=\bigg| \sum_{|m|\geq M+1} H_m^{(1)}(k|x-x_j|)\text{exp}[\rmi m\theta]k\frac{(y-x_j)\cdot \nu}{|y-x_j|^3}J_m'(k|y-x_j|)
     \\&+H_m^{(1)}(k|x-x_j|)\text{exp}[\rmi m\theta]J_m(k|y-x_j|)(-\rmi m) \frac{\partial}{\partial \nu}\text{arg}(y-x_j)\bigg|\\
       &\leq  \sum_{|m|\geq M+1} |H_m^{(1)}(k|x-x_j|)k\frac{(y-x_j)\cdot \nu}{|y-x_j|}J_m'(k|y-x_j|)|
     \\&+|H_m^{(1)}(k|x-x_j|)J_m(k|y-x_j|)m\frac{\partial}{\partial \nu}\text{arg}(y-x_j)|.
\end{align*}
Since $H^{(1)}_{-m}=(-1)^mH^{(1)}_m$ and $J_{-m}=(-1)^mJ_m$ we can reduce the sum to
\begin{align}
    R'_{j,M}(x;k)&\leq 2\sum_{m=M+1}^\infty |H_m^{(1)}(k|x-x_j|)k\frac{(y-x_j)\cdot \nu}{|y-x_j|}J_m'(k|y-x_j|)| \nonumber
     \\&+2\sum_{m=M+1}^\infty |H_m^{(1)}(k|x-x_j|)J_m(k|y-x_j|)m\frac{\partial}{\partial \nu}\text{arg}(y-x_j)|.\label{eq.sumpsplit}
\end{align}
We deal with the two latter sums separately.
We start with the second sum due to its similarity to the monopole error. Noting that 
\begin{equation}\label{eq.domarg}
\bigg| \frac{\partial}{\partial \nu}\text{arg}(y-x_j)\bigg| = \bigg|\frac{(y-x_j)_\perp}{|y-x_j|^2}\cdot \nu \bigg| \leq \frac{1}{|y-x_j|} \leq C  
\end{equation}
where  for a vector $u=(u_1,u_2)\in \mathbb{R}^2$, $u_\perp=(-u_2,u_1) $ and the positive constant $C$ is defined by $C=\max_{y\in Y}|y-x_j|^{-1}$.
Thus, combining \eqref{eq.Jnbis-rem}, \eqref{eq.Hn-rem} and \eqref{eq.domarg} gives that 
\begin{equation}\label{eq.domsecondsumbis}
\sum_{m=M+1}^\infty m |H_m^{(1)}(k|x-x_j|)J_m(k|y-x_j|)\frac{\partial}{\partial \nu}\text{arg}(y-x_j)|\leq C_3 \sum_{m=M+1}^\infty a_x^m,
\end{equation}
where the positive constant $C_3$ depends only on $Y$, $X$ and $K$.\\
Now, we estimate the first sum of \eqref{eq.sumpsplit}. By virtue of the estimates \eqref{eq.Jnprimebis-rem} and \eqref{eq.Hn-rem}, one gets that there exists a constant $C_4>0$ such that
\begin{equation*}
\sum_{m=M+1}^\infty |H_m^{(1)}(k|x-x_j|)k\frac{(y-x_j)\cdot \nu}{|y-x_j|}J_m'(k|y-x_j|)| \leq \frac{C_4}{|x-x_j|}  \sum_{m=M+1}^\infty a_x^m.
\end{equation*}
Setting  $C_5=C_4 \, (\max_{x\in X}  |x-x_j|)^{-1}>0$, one obtains 
\begin{equation}\label{eq.domfirstsum}
\sum_{m=M+1}^\infty |H_m^{(1)}(k|x-x_j|)k\frac{(y-x_j)\cdot \nu}{|y-x_j|}J_m'(k|y-x_j|)| \leq  C_5\sum_{m=M+1}^\infty a_x^m.
\end{equation}
Combining \eqref{eq.sumpsplit}, \eqref{eq.domsecondsumbis} and \eqref{eq.domfirstsum} yields  the second inequality of \eqref{eq.graf}.
\end{proof}

\subsection{Numerical experiments for truncation error estimates}
\label{sec:ten}

We illustrate our bounds numerically in the case where there is only one source to move, i.e. $Y = \{ y \}$. The bounds in \cref{pro:grafs_trunc} involve a quantity $a_x$ that can be estimated from the relative positions of $x$, $y$ and $x_j$ (repectively, the evaluation point and the original and new source positions).  The bounds also involve constants $C_1, C_2$ that may depend in non-obvious ways on the different choices of compact sets in space and wavenumber. To estimate $C_1$ (resp. $C_2$) for a particular choice of $(x,k) \in X\times K$,  we assume the truncation error has the form predicted by the respective upper bound in \eqref{eq.graf}, and we find the $C_1$ (resp. $C_2$) that matches the actual error explicitly for one small value of $M$.  We repeat this estimate on a grid for $X \times K$ and then take the maxima of the estimates for $C_1$ (resp. $C_2$)  over the grid.

\begin{figure}
\centering
 \begin{tabular}{c}
    \includegraphics[width=55mm]{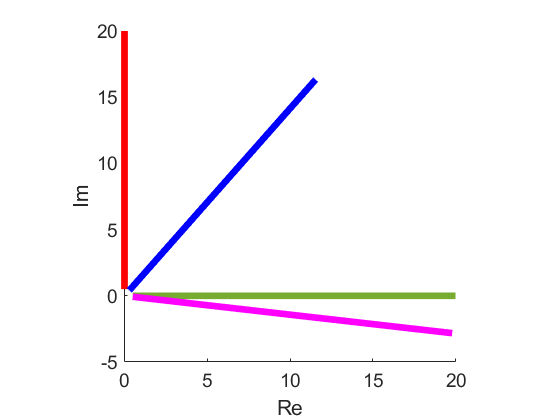}\\
    \begin{tabular}{l|l|l}
    Condition & Range & Color\\\hline
    real  & $\mathcal{K}_1 = [0.5,20]$ & green\\
    imaginary & $\mathcal{K}_2 = \{\rmi\theta ~|~  \theta \in \mathcal{K}_1 \}$ & red\\
    dissipative & $\mathcal{K}_3 = \{(1+\sqrt{2}\rmi)\theta/\sqrt{3} ~|~ \theta \in \mathcal{K}_1  \}$ & blue\\
    amplifying & $\mathcal{K}_4 = \{(99-\rmi\sqrt{199}\theta)/100 ~|~ \theta \in \mathcal{K}_1 \}$ & pink
    \end{tabular}\\
\end{tabular}
    \caption{Wavenumber ranges used in the numerical experiments  and their visualization in the complex plane.}
    \label{fig:cases}
\end{figure}

\begin{figure}
    \centering
    \begin{tabular}{cc}
    \includegraphics[width=55mm]{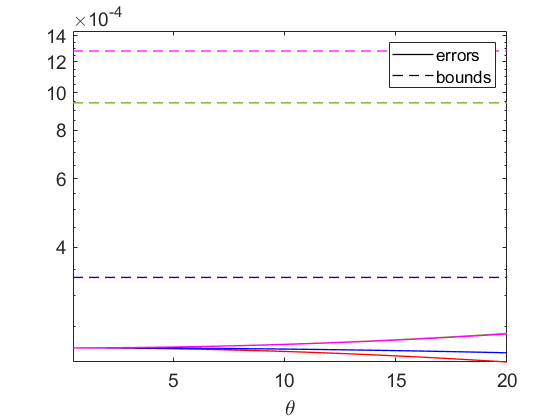} & \includegraphics[width=55mm]{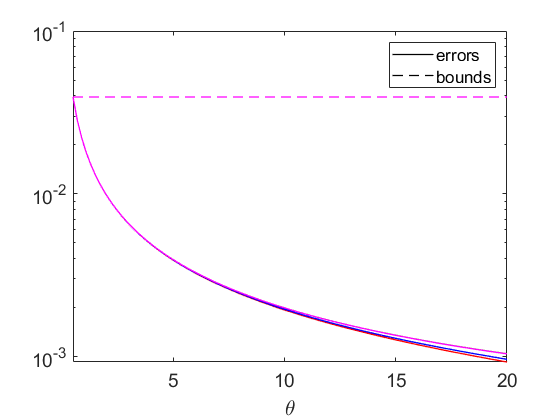}\\
    \captionfont (a) Monopoles & \captionfont (b) Dipoles\\
    \end{tabular}
    \caption{The monopole (a) and dipole (b) errors (logarithmic scale), $R_{j,M}$ and $R'_{j,M}$ in \eqref{eq:rjm} respectively, and the bounds from \cref{pro:grafs_trunc} for a single point in space which was relatively close to the source. The color represents the set for $k$ as given in \cref{fig:cases}. For the dipoles we cannot differentiate the different bounds as they all lie on top of each other. We observe that the dipole errors and bounds are larger. This is in line with  \cref{pro:grafs_trunc}: the dipole error bound decays slower than the monopole error bound.}
    \label{fig:freq_err}
\end{figure}

\begin{figure}
    \centering
    \begin{tabular}{cc}
    \includegraphics[width=55mm]{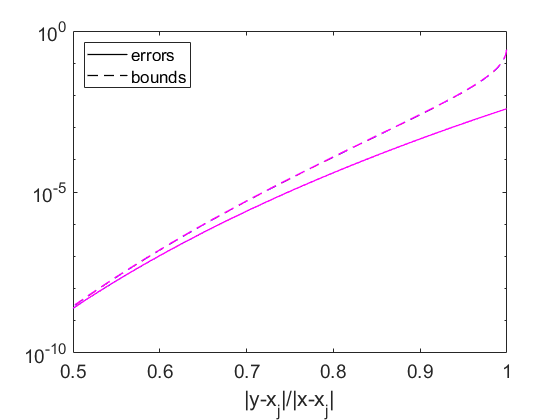} & \includegraphics[width=55mm]{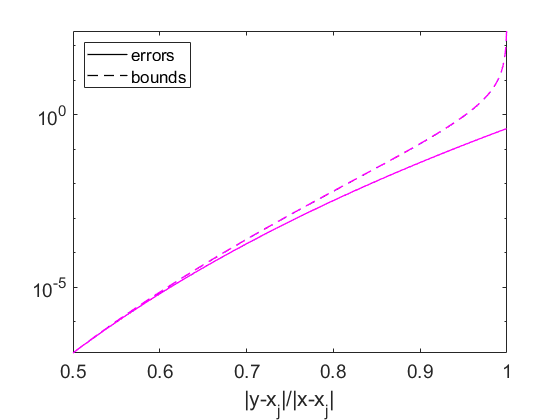}\\
    \captionfont (a) Monopoles &\captionfont  (b) Dipoles 
    \end{tabular}
    \caption{The monopole (a) and dipole (b) errors (logarithmic scale), $R_{j,M}$ and $R'_{j,M}$ in \eqref{eq:rjm} respectively, and bounds for a range of points in space at a real $k=1$, pure imaginary $k=\rmi$, complex (dissipative) $k=1/\sqrt{3}+\rmi\sqrt{2}/\sqrt{3}$, and complex (gain) $k=(99-\rmi\sqrt{199})/100$ wavenumbers with the colors represented in \cref{fig:cases}. The plots for different wavenumbers cannot be differentiated.}
    \label{fig:space_err}
\end{figure}

In \cref{fig:freq_err} we show these bounds for $\{x\} \times K$, for different choices of wavenumber sets in complex plane and for the fixed evaluation point $x = (0,0.43)$. The original source location is $y= (0,0)$ and it is moved to the new location $x_j =(0,0.2)$. We took $M=4$ to approximate the constants $C_1$ and $C_2$ over $\{x\} \times K$ and then use our estimated $C_1$ and $C_2$ to predict the truncations errors with $M=20$ terms. The different wavenumber ranges in the complex plane that we considered are summarized in  \cref{fig:cases}.
Then in \cref{fig:space_err}, we estimated the constants $C_1$ and $C_2$ on $X \times \{k\}$ for four different wavenumbers $k \in \complex$. Here $X$ is the region $ X = \{ x \in \real^2 ~|~ 1/2 \leq |y-x_j| / |x-x_j| \leq 1 \}$, i.e. the annulus for which the ratio in the geometric series ansatz belongs to $[1/2,1]$.

\section{Active exterior cloaking at fixed frequency}
\label{sec:aec}
One can achieve active cloaking by observing  \cite{Miller:2005:OPC} that a distribution of monopoles and dipoles on the boundary $\partial \Omega$ of a bounded open $\Omega \subset \real^2$ can create a field $u_c$ that cancels out the incident or probing field $u_i$ inside a region $\Omega$, while vanishing outside, or in  other words satisfying \eqref{eq:ucfield}. By applying the Green identities (see e.g. \cite{Colton:2013:IAE}) or potential theory (see e.g. \cite{Kellogg:1967:FPT}) the function $u_c=-u_i$ is given for $x \notin \partial \Omega$ by
\begin{equation}\label{eqn:brepfreq}
    u_c(x;k) = \int_{\partial \Omega} dS(y) [\frac{-\partial u_i}{\partial \nu(y)}(y;k)G(x-y;k)+u_i(y;k)\frac{\partial G}{\partial \nu(y)}(x-y;k)],
\end{equation}
where $G(x;k)$ is the Green function \eqref{eq:greenfun}.
\begin{Rem}
\label{rem:validity}
The representation formula \eqref{eqn:brepfreq} is valid for example when $\Omega$ has Lipschitz boundary $\partial\Omega$. To see this, we assume $\Omega \subset \mathcal{O}$ where $\mathcal{O}$ is an open set and the incident field $u_i\in H^1_{loc}(\mathcal{O})$ solves (in the distributional sense) the homogeneous Helmholtz equation $\Delta u_i+k^2 u_i=0$ in $\mathcal{O}$ for $k\in \mathbb{C}\setminus \{0\}$. 
Then as $-\Delta u_i=k^2 u_i$ on $\mathcal{O}$, one easily proves by interior elliptic regularity of the minus Laplacian  operator (applying iteratively  Theorem 2 page 314 of \cite{Evans:2010:PDE}) that $u_i\in C^{\infty}(\overline{\Omega})$. 

We point out  that $u_c=-u_i\in C^{\infty}(\overline{\Omega})$ and the outward normal vector $\nu(y)\in L^{\infty}(\partial \Omega)$ since $\partial \Omega$ is a Lipschitz boundary. Thus, it is clear that the Dirichlet trace $u_i$ is smooth on $\partial \Omega$ and that the Neumann trace $[\partial u_i/\partial \nu(y)](y;k)$ is in $L^\infty(\partial \Omega)$. Hence,  the integrand in \eqref{eqn:brepfreq} is integrable as a sum of two products of $L^\infty(\partial\Omega)$ functions. Indeed since we have $x\notin \partial \Omega$,  the Green function $G(x-y;k)$ is smooth for $y \in \partial \Omega$  and its normal derivative $[\partial G/ \partial \nu(y)](x-y;k)$ is in $L^\infty(\partial \Omega)$ as a function of $y$.

\end{Rem}
 
To get \emph{exterior} cloaking, the idea is to move the monopoles and dipoles on the portions $\partial\Omega_j$ of the boundary $\partial\Omega$ to the new source locations $x_j$. Formally, this can be done by replacing the Green function and its normal derivative in the representation formula \eqref{eqn:brepfreq} by their series expansions \eqref{eqn:grafs} and \eqref{eqn:grafs.dipole}. \Cref{thm:conv} and \cref{rem:validity} allow to permute the order of the series and the integral over $\partial\Omega$ (since the series is normally convergent with respect to $y$). Thus we can express the new cloaking field as 
\begin{equation}
    u_e(x;k) = \sum_{j = 1}^{N_{dev}} \sum_{m= -\infty}^\infty b_{j,m}V_m(x-x_j;k),
    \label{eq:ue}
\end{equation}
where $V_m$ are multipolar sources \eqref{eq:vm} and the coefficients $b_{j,m}$ are given by \eqref{eq:bjm} in terms of integrals over the $\partial\Omega_j$, identical to those obtained in \cite{Guevara:2009:AEC}. We emphasize that \cref{thm:conv} is valid for complex $k$ with the exception of the negative real axis, whereas the result in \cite{Guevara:2009:AEC} is only proven for real $k$ positive. Moreover, \cref{thm:conv} leverages on the Graf addition formula truncation error estimates in \cref{pro:grafs_trunc}, to give the truncation error when we consider instead the truncated fields:
\begin{equation}\label{eqn:ext_trunc}
    u^{(M)}_e(x;k) = \sum_{j = 1}^{N_{dev}} \sum_{m= -M}^M b_{j,m}V_m(x-x_j;k).
\end{equation}
This error estimate is novel and applies to the results in \cite{Guevara:2009:AEC}.

\begin{theorem}\label{thm:conv} Let $\Omega\subset \real^2$ be a bounded open set with Lipschitz boundary $\partial\Omega$. Assume $u_i$ is a $H^1_{loc}(\mathcal{O})$ solution to the homogeneous Helmholtz equation, where $\mathcal{O}$ is an open set containing $\overline{\Omega}$.
Define the region $R = R_1 \cup \cdots \cup R_{N_{dev}}$, i.e. the union of the disks $R_j$ in \eqref{eq:rj}. Let $K$ be a compact subset of $\mathbb{C} \setminus (-\infty,0]$ and $X$ a compact subset of $\real^2 \setminus R$. Define the coefficients $b_{j,m}$ in \eqref{eq:ue} and \eqref{eqn:ext_trunc} by
\begin{equation} 
b_{j,m} = \int_{\partial \Omega_j} dS(y)\bigg[-\frac{\partial u_i}{\partial \nu(y)}(y;k)U_m(y-x_j;k) + u_i(y;k)\frac{\partial U_m(y-x_j;k)}{\partial \nu(y)}\bigg],
\label{eq:bjm}
\end{equation}
where $j=1,..,N_{dev}$, $m \in \mathbb{Z}$, and
$U_m(x;k) = 
    J_m(k|x|)\exp[-\rmi m\text{arg}(x)]$.
Then there exists a constant $C>0$ (which may depend on $K$, $X$, $u_i$ and $\partial u_i / \partial \nu(y)$) such that for any $(x,k) \in X \times K$, 
\[\left|u_c(x;k)-u^{(M)}_e(x;k)\right| \leq C \frac{a^{M+1}}{1-a}\]
where
\[
 a = \max_{x\in X} \max_{j=1,\ldots,N_{dev}}  \max_{y \in \partial\Omega_j} \frac{|y - x_j|}{|x - x_j|}<1.
\]
In particular for any $x \notin R$ and $k\notin (-\infty,0]$, we have $u_c(x;k) = u_e(x;k)$.
\end{theorem}

\begin{Rem}
 The integral appearing in the definition of the $b_{j,m}$ in \eqref{eq:bjm} can be expressed as a series when the incident field $u_i$ is given in terms of its cylindrical wave expansion \cite[Theorem 2]{Norris:2012:SAA}. Although the series expansion is proven for the 2D Helmholtz equation with $k>0$, we conjecture it is valid for complex $k$.
\end{Rem}

\begin{Rem}
 Controlling fields outside of a bounded open set $\Omega$ can be useful for the mimicking problem (making a scatterer inside $\Omega$ look like another one) or for cloaking a source inside $\Omega$. This requires an \emph{exterior} version of the Green representation formula \eqref{eqn:brepfreq}, which is valid for $\Im(k) \geq 0$, when the field to reproduce $u_i$ is a solution to the homogeneous Helmholtz equation outside of $\Omega$ and satisfies the Sommerfeld radiation condition (see e.g. \cite[Theorem 3.3]{Colton:1983:IEM}). We  are not aware of the validity of this result for gain media ($\Im (k) <0$). Therefore we anticipate that \cref{thm:conv} can be adapted to control fields outside of $R$ for $\Im(k) \geq 0$ and $k \notin (-\infty,0]$.
\end{Rem}

\subsection{Proof of \cref{thm:conv}}
\begin{proof}
We rewrite the boundary integral representation \eqref{eqn:brepfreq} of the cloaking field $u_c$ as a sum of integrals over the portions $\partial\Omega_j$ of the boundary. We then apply \eqref{eqn:grafs} to yield 

\begin{equation}
\begin{split}
    u_c(x)&= \sum_{j=1}^{N_{dev}} \int_{\partial \Omega_j} dS(y) \bigg[\frac{-\partial u_i}{\partial \nu(y)}(y)\frac{\rmi}{4}\sum_{m=-\infty}^\infty V_{m}(x-x_j)U_m(y-x_j)\\
    &+u_i(y)\frac{\rmi}{4}\frac{\partial}{\partial \nu(y)}\sum_{m=-\infty}^\infty V_{m}(x-x_j) U_m(y-x_j) \bigg]
\end{split}
\label{eq:ucexpanded}
\end{equation}
which holds for $x \notin R$. We approximate the cloak field by $u_e^{(M)}$ with coefficients $b_{j,m}$ chosen as in \eqref{eq:bjm} to match the $|m|\leq M$ terms in the series in \eqref{eq:ucexpanded}. Thus the error we make by approximating $u_c(x)$ by $u_e^{(M)}$ at some $x \notin R$ can be bounded by
\[\begin{split} |u_c(x) &- u_e^{(M)}(x)| =\\ &\bigg|\sum_{j=1}^{N_{dev}} \int_{\partial \Omega_j} dS(y) \bigg[\frac{-\partial u_i}{\partial \nu(y)}(y)\frac{\rmi}{4} \sum_{|m| \geq M+1} V_{m}(x-x_j)U_m(y-x_j)\\
    &+u_i(y)\frac{\rmi}{4}\frac{\partial}{\partial \nu(y)}\sum_{|m| \geq M+1} V_{m}(x-x_j) U_m(y-x_j) \bigg] \bigg|\\
    &\leq \sum_{j=1}^{N_{dev}} \int_{\partial \Omega_j} dS(y) \bigg[ \bigg|\frac{-\partial u_i}{\partial \nu(y)}(y)\frac{\rmi}{4} \bigg| \bigg| \sum_{|m| \geq M+1} V_{m}(x-x_j)U_m(y-x_j) \bigg|\\
    &+\bigg|u_i(y)\frac{\rmi}{4} \bigg| \bigg|\frac{\partial}{\partial \nu(y)}\sum_{|m| \geq M+1} V_{m}(x-x_j) U_m(y-x_j)\bigg| \bigg].
\end{split} \]
We notice that 
\[ 
\begin{split}
    R_{j,M}(x;k) &=\bigg| \sum_{|m| \geq M+1} V_{m}(x-x_j)U_m(y-x_j) \bigg|\\
    R'_{j,M}(x;k) &= \bigg|\frac{\partial}{\partial \nu(y)}\sum_{|m| \geq M+1} V_{m}(x-x_j) U_m(y-x_j)\bigg|,
\end{split}
\]
allowing us to apply \cref{pro:grafs_trunc} to bound the truncation by remainders of a geometric series. It follows that there exists two positive constants $C_{1,j}$ and $C_{2,j}$ that depend on the compact sets $X$, $Y_j=\partial \Omega_j$ and $K$ such that  
\begin{equation*}
    |u_c(x) - u_e^{(M)}(x)| \leq  \sum_{j=1}^{N_{dev}} \int_{\partial \Omega_j} dS(y)\bigg( \bigg|\frac{\partial u_i}{\partial \nu(y)}(y)\frac{1}{4}C_{1,j}\bigg| + \bigg|u_i(y)\frac{1}{4}C_{2,j} \bigg|\bigg)\frac{a_{x,j}^{M+1}}{1-a_{x,j}},
\end{equation*}   
with 
$$
a_{x,j}=  \frac{\max_{y\in \partial_{\Omega_j}} |y-x_j|}{|x-x_j|}.
$$   
Setting $C_1=\displaystyle\max_{j=1, \ldots, N_{dev}} C_{1,j} $ and $C_2=\displaystyle\max_{j=1, \ldots, N_{dev}} C_{2,j} $, it follows that 
\begin{equation*}
    |u_c(x) - u_e^{(M)}(x)|\leq  n \max_{j}|\partial \Omega_j| \, \tilde{C} \, \frac{a^{M+1}}{1-a},
    \label{eq:error:bound}
  \end{equation*}  
where from \cref{rem:validity} we have that $u_i \in C^0(\partial\Omega)$ and $\partial u_i / \partial \nu(y)$ is $L^\infty(\partial\Omega)$:
\begin{align}
    \tilde{C} &= \sup_{y \in \partial \Omega}\bigg|\frac{\partial u_i}{\partial \nu(y)}(y)\frac{1}{4}C_1\bigg|+ \max_{y \in \partial \Omega}\bigg|u_i(y)\frac{1}{4}C_2 \bigg|,~\text{and}\\
    a &= \max_{x\in X} \max_{j=1,\ldots,n}  a_{x,j}<1,
    \label{eq:a}
\end{align}
and the error bound follows by letting $C=N_{dev}\tilde{C}\max_j|\partial \Omega_j|.$  In addition, we have for $x \notin R$ and $k \in \complex \setminus (\infty,0]$ that  $|u_c(x)-u_e(x)|=\lim_{M\to \infty}|u_c(x)-u_e^{(M)}(x)|=0$ since $a<1$. Hence we have $u_c=u_e$ for $x \notin R$. The fields do not agree for $x \in R$, because for at least one $j \in \{1,\ldots,N_{dev}\}$, the series in \eqref{eq:ue} diverges. 
\end{proof}

\subsection{Stability through the maximum principle}
\label{sec:stability}
Homogeneous Helmholtz equation solutions satisfy a strong maximum principle if $\Re(k^2) < 0$ or equivalently
\begin{equation}\label{def.strongdissipation}
|\Im(k)|>|\Re(k)|. 
\end{equation}
Although we use this result for constant isotropic media, it has been proved in the very general context of the Helmholtz equation with anisotropic heterogeneous media \cite[corollary 2.1]{Kresin:1993}. Another proof in the case of isotropic heterogeneous media appears in  \cite[theorem 6]{Krutitskii:1998:TDP}.

We now state the strong maximum principle. By interior regularity (see \cref{rem:validity}), a homogeneous solution $u \in H^1_{loc}(\mathcal{O})$ to the Helmholtz equation  in an open set $\mathcal{O}$  with a wavenumber $k$ satisfying \eqref{def.strongdissipation}, is $C^\infty(\overline{B})$ on any open bounded subset $\mathcal{B}$ with Lipschitz boundary $\partial B$ satisfying $\overline{\mathcal{B}} \subset \mathcal{O}$.
Thus, one can apply the strong maximum principle on the set $\mathcal{B}$ to get on one hand that
$$
\max_{x \in \overline{\mathcal{B}}} |u(x)|=\max_{x \in \partial \mathcal{B}} |u(x)|,
$$
and on the other hand that the maximum of $|u|$ is only reached on the boundary $\partial \mathcal{B}$ of $\mathcal{B}$. We note that the \emph{strong} maximum principle is not valid for $k$ outside the region \eqref{def.strongdissipation} as one can find examples of solutions violating it \cite{Kresin:1993}.

In particular if $u$ and $v$ are smooth solutions to the homogeneous Helmholtz equation in $\mathcal{B}$ with $k$ satisfying \eqref{def.strongdissipation}, the maximum of the error $|u(x) - v(x)|$ is attained \emph{only} at the boundary $\partial \mathcal{B}$. In other words, the error within the domain is controlled by the error on the boundary (the Dirichlet data). This can be viewed as a form of stability for the boundary integral representation \eqref{eqn:brepfreq}. Moreover, $u_e$ and $u_e^{(M)}$ are $C^\infty(\overline{\mathcal{B}})$ solutions  to the homogeneous Helmholtz equation on $\mathcal{B}$, where $\mathcal{B}$ is a bounded open set such that $\overline{\mathcal{B}} \subset \real^2 \setminus R$ (see \cref{thm:conv}). Therefore we can conclude from the maximum principle that the truncation error of exterior cloaking ($|u_e(x) - u_e^{(M)}(x)|$) reaches its maximum over $\overline{\mathcal{B}}$ only on the boundary $\partial\mathcal{B}$. 
Finally we point out that when numerically evaluating the boundary representation formula \eqref{eqn:brepfreq}, we use finitely many monopole and dipole sources on the domain $\partial \Omega$. Following the same argument, the error we make with this discretization is also maximum  on the boundary of any bounded domain $\mathcal{B}$ such that $\overline{\mathcal{B}} \subset \real^2 \setminus \partial \Omega$. We numerically illustrate in \cref{fig:max} that the maximum principle predicts that the maximum cloaking errors occurs on the boundary of a region and not inside.

\subsection{Numerical experiments}\label{sec:numericsfreq}
We explain how we evaluate the truncated cloak field $u_e^{(M)}$ in \cref{sec:evaluation}. Then the truncation errors are predicted in \cref{sec:errorbound} using the error bounds in \cref{thm:conv}. Finally we explain in \cref{sec:scattered} how we calculate scattered fields when $\Im(k) \geq 0$.  

\subsubsection{Evaluation of the cloak field}

\begin{figure}
    \centering
    \begin{tabular}{cc}
    \includegraphics[width=0.35\textwidth]{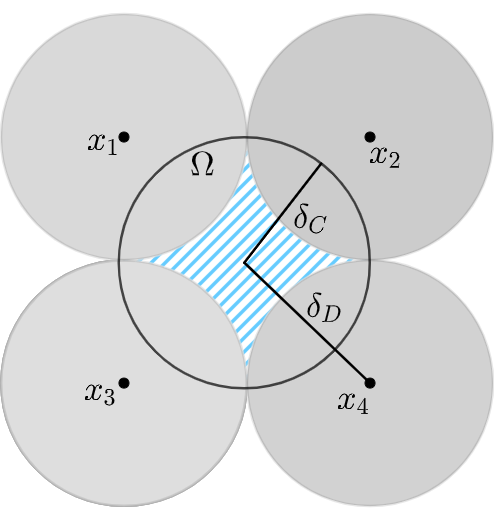} & \includegraphics[width=0.35\textwidth]{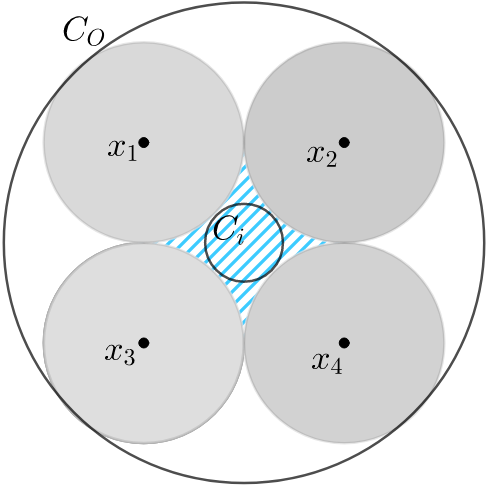}\\
    \captionfont(a) Inner and exterior cloak configuration &\captionfont(b) Cloak with in/circumscribed circles
    \end{tabular}
    \caption{(a) The configuration of the exterior cloak used for four exterior sources ($x_j$) to maximize the region where an object can be dissimulated (the blue striped region), we take $\delta_D = \sqrt{2}\delta_C.$ (b) Cloaking region with smallest inscribed circle with radius $r_{C_i}$ and largest circumscribed circle with radius $r_{C_O}.$}
    \label{fig:cart_config}
\end{figure}

\begin{figure}
    \centering
    \begin{tabular}{ccc}
    & Real part of field & Real part of reproduction\\
    \raisebox{4em}{\rotatebox{90}{$k=10$}} &\includegraphics[width=55mm]{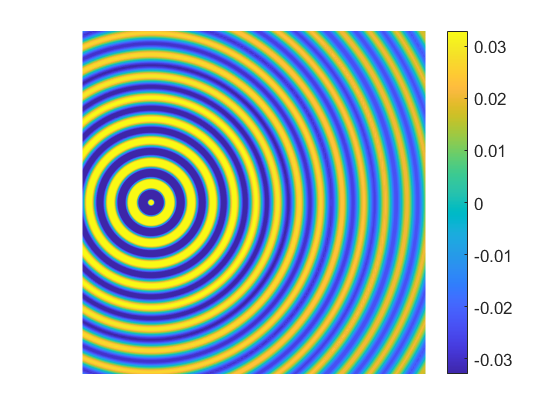} & \includegraphics[width=55mm]{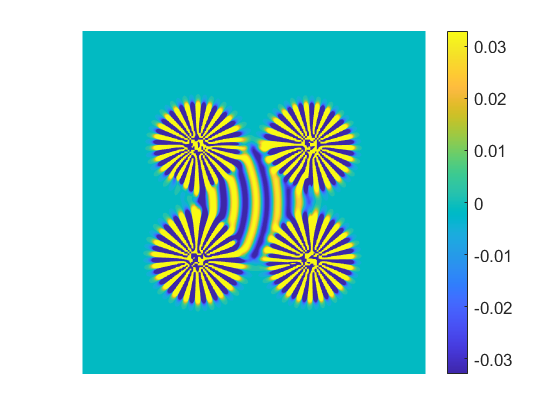}\\
     \raisebox{4em}{\rotatebox{90}{$k=\rmi/2$}} &\includegraphics[width=55mm]{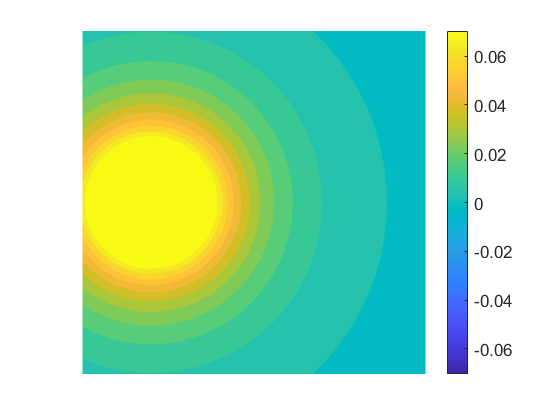} & \includegraphics[width=55mm]{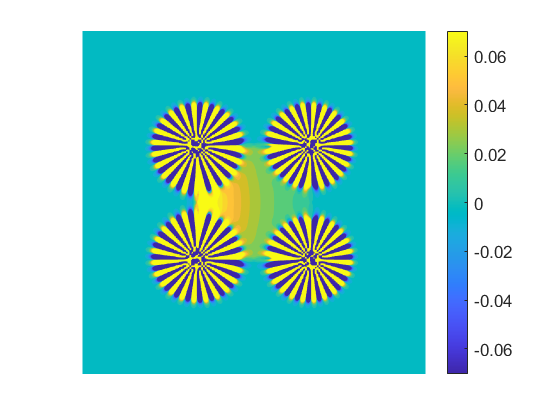}\\
     \raisebox{4em}{\rotatebox{90}{$k=10+\rmi/2$}} &\includegraphics[width=55mm]{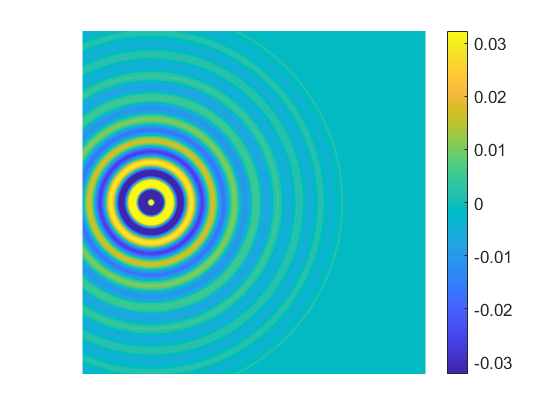} & \includegraphics[width=55mm]{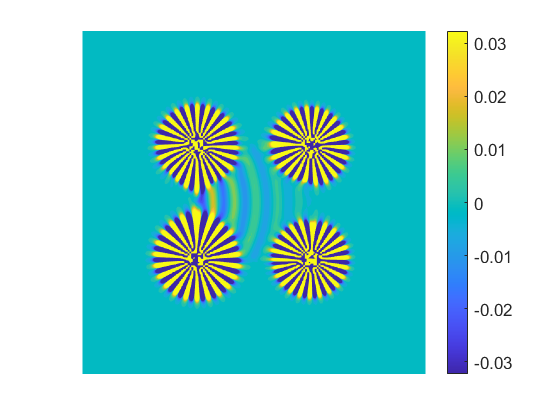}\\
     \raisebox{4em}{\rotatebox{90}{$k=10-\rmi/2$}} &\includegraphics[width=55mm]{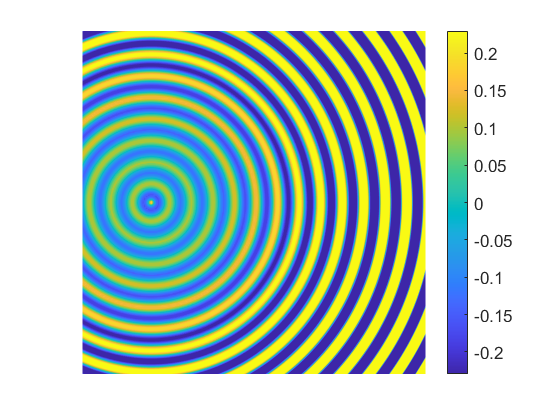} & \includegraphics[width=55mm]{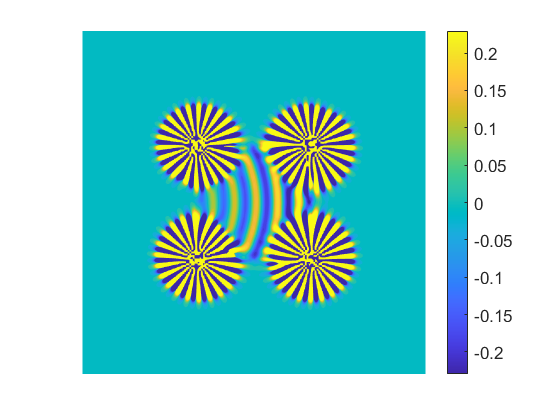}
    \end{tabular}
    \caption{Field reproductions (right) and the original field (left) at different wavenumbers on the square $[0,1]^2$ with a point source located at $(2,5)$. The color scale was kept the same for each $k$ and was chosen to highlight the different behaviors of the point sources for different $k$.}
    \label{fig:rep}
\end{figure}

\label{sec:evaluation}
We illustrate \cref{thm:conv} numerically using a disk region $\Omega$ and $n=4$ sources, as shown in \cref{fig:cart_config}.  While we chose to illustrate exterior cloaking  with four multipolar sources, only three are necessary in two dimensions to give a non-empty region cloaking \cite{Guevara:2011:ECA}.  Cloaking fields $u_e^{(M)}$ with $M=22$ are shown in \cref{fig:rep} for several representative wavenumbers on the square $[0,10]^2$ using a $200 \times 200$ uniform grid. The disk $\Omega$ is centered at $(5,5)$ and with radius $\delta_C = 10/6$. The $x_j$ are placed on a circle of radius $\delta_D$ (see \cref{fig:cart_config}), where $\delta_D$ is chosen to maximize the area of the cloaking region for $\delta_C$ fixed. The optimization is done via a simple geometric argument similar to \cite{Guevara:2011:ECA} and gives $\delta_D = 5\sqrt{2}/3$. The incident field is generated by a point source at $y=(2,5)$.
To evaluate the truncated exterior cloaking field, $u_e^{(M)}$, we use an equi-spaced discretization of $\partial \Omega$, into points $y_i$ with $i=1,...,n_{int}$. We split $\partial \Omega$ into $n=4$ regions each associated with  a new source location $x_j$.  We choose $n_{int}$ so that there is an equal number of discretization points of $\partial\Omega$ for each $\partial\Omega_j$ and choose the $x_j$ such that $\max_{y_i \in \partial \Omega_j}|y_i-x_j|$ is equal for all $j$ in order to keep the size of the theoretical divergence regions $R_j$ of our devices equal. The integrals over $\partial \Omega_j$ that determine the coefficients $b_{j,m}$ in \cref{thm:conv} are approximated using the midpoint rule (so that the total integral over $\partial\Omega$ is the trapezoidal rule).

We note that the color scale in \cref{fig:rep} is deliberately limited to exclude the large fields near the new source locations $x_j$ which are due to the singularity of $u_e^{(M)}$ at the $x_j$. This may seem as an impediment to physically realize such cloaking devices. However, as noted in \cite{Guevara:2011:ECA}, it is possible to use the Green exterior representation formula (valid for $\Im(k)\geq 0$, see \cite[Theorem 3.3]{Colton:1983:IEM}) to replace the multipolar sources by a distribution of monopoles and dipoles on some boundary enclosing each of the $x_j$. Since the cloak field $u_e^{(M)}$ is smaller, we expect it is easier to realize in practice. The drawback is that these ``extended cloaking devices'' leave only small gaps between the cloaked region and the exterior. By \cref{thm:conv} we expect that $|u_e^{(M)}(x)|\to \infty$ for $x \in R$, so increasing $M$ leads to smaller gaps. So there is a tradeoff between getting larger gaps and approximating the ideal cloak field $u_e$.

\subsubsection{Computation of error bounds and the maximum principle}
\label{sec:errorbound}

\begin{figure}
	\centering
	\includegraphics[width=55mm]{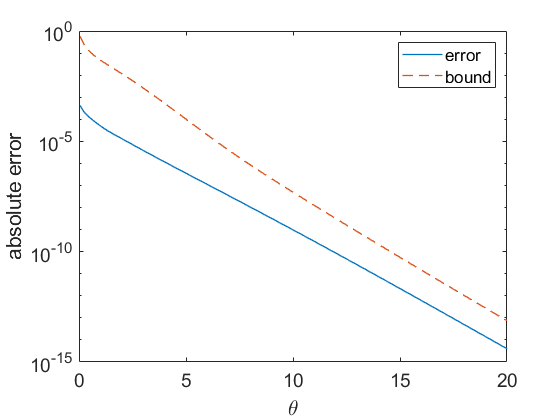}
	\caption{Predicted and actual maximum error (logarithmic scale) on circles of radius slightly larger than $r_{c_O}$ and smaller than $r_{c_i}$, which are  outside of the divergence region of the cloak field. The abscissa corresponds to the parameter $\theta$ for the segment of wavenumbers $\mathcal{K}_3$, as defined in \cref{fig:cases}.  }
	\label{fig:freqvtrunc}
\end{figure}

In order to use the error bounds from \cref{thm:conv} in our particular geometric setup, we define the radius $r_{C_i}$ (resp. $r_{C_O}$) of the largest (resp. smallest) inscribed (resp. circumscribed) circle that is inside (resp. outside) the divergence region $R$ (defined in \cref{thm:conv}). The inscribing and circumscribing circles are represented in \cref{fig:cart_config}. We recall from \eqref{eq:error:bound} that the cloak field truncation error can be bounded by the truncation error of a geometric series with ratio $a$ that is determined by the relative positions of the $\partial \Omega_j$, the $x_j$ and the region of interest where we want to evaluate the fields, see \eqref{eq:a}. Since we expect the cloaking fields to diverge close to $R$, it does not make sense to evaluate the errors on the inscribing and circumscribing circles. We do it instead on slightly smaller or larger circles of radii $r_{C_0} + 0.1\delta_C$ and $r_{C_i}-0.1\delta_C$. If we take the region $X$ from \cref{thm:conv} to be the union of these two circles, then a simple geometric argument yields that there are eight points in $X$ that attain the maximum over $X$ in the definition of $a$ \eqref{eq:a}. At each of these points the ratio of the geometric series ansatz is the same, so we can conclude the truncation error can be bounded by $C (1-a^{M+1})/(1-a)$, where $a<1$, but the constants $C$ depend on the point. We first estimate the constant $C$ at a point $x$ by using the ``empirical method'' we used in \cref{sec:ten}. In other words we find the $C$ for which $|u_e^{(3)}(x) - u_e(x)|$  is equal to $C (1-a^{M+1})/(1-a)$ with $M=3$. Then we take the worst case scenario, i.e. the largest of such $C$ for the eight points in $X$ that we considered. We emphasize that this is a heuristic meant to simplify the exhaustive method, where would have to evaluate the largest $C$ for all $x\in X$. We summarize in \cref{fig:freqvtrunc} the application of this heuristic for wavenumbers $k\in \mathcal{K}_3$ (as defined in \cref{fig:cases}). In these experiments we used $128$ equispaced discretization points for $\partial\Omega$,  $\delta_C= 10/6$ and $\delta_D = 5\sqrt{2}/3$. Finally the incident field  we used for this experiment was a point source located at $x=(8,5)$. As can be seen from \cref{fig:freqvtrunc}, the error bound we obtain for $M=22$ overestimates the actual error and follows the same trend for varying wavenumber.

\begin{figure}
 \centering
 \begin{tabular}{cc}
\includegraphics[width=60mm]{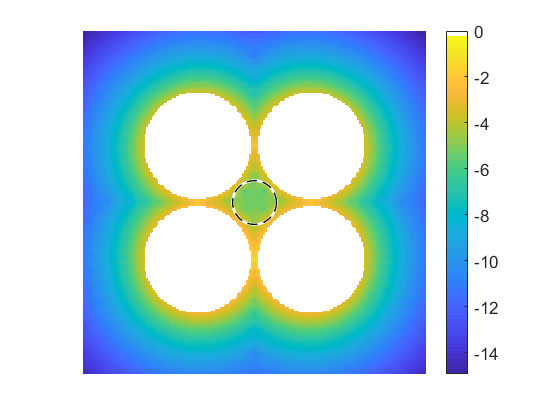} & \includegraphics[width=60mm]{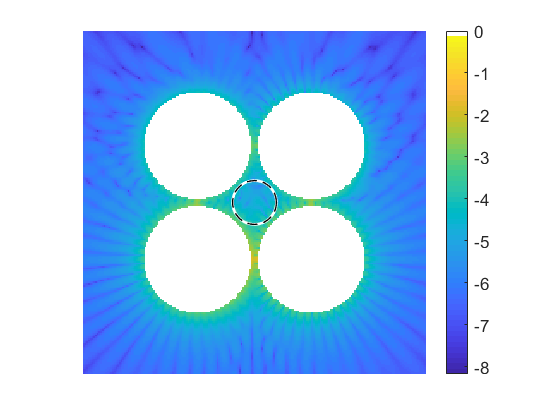}\\
 (a) & (b)
 \end{tabular}
 \caption{We display the cloaking field truncation error $\log_{10} |u_e^{M} - u_e|$ (outside of the ``extended cloaking devices'' in white) corresponding to (a) $k=\rmi/2$ and (b) $k= 10+\rmi/2$. We note that (a) (resp. (b)) corresponds to the second (resp. third) row in \cref{fig:rep}. By applying the maximum principle on the disk $\mathcal{B}$ (dashed curve), we see that the maximum error is attained on $\partial\mathcal{B}$ in (a) but not in (b). The difference is that the wavenumber in (a) satisfies $|\Im(k)| > |\Re(k)|$ so a version of the maximum principle applies, see also \cref{fig:kdiagram}.}
 \label{fig:max}
\end{figure}

We illustrate in \cref{fig:max} that when $|\Im(k)|>|\Re(k)|$, the maximum principle (\cref{sec:stability}) can be used to predict where the maximum cloaking error occurs. In fact the wavenumbers we used in for \cref{fig:freqvtrunc} also allow us to use the maximum principle to observe that a bound for the truncation error on the boundary of the circle with radius $r_{C_i} - 0.1\delta_C$ automatically leads to a bound on the whole disk of same radius.

\subsubsection{Calculating scattered fields for $\Im(k) \geq 0$}
\label{sec:scattered}
To demonstrate cloaking, we recall how to calculate scattered fields from a sound-soft (or homogeneous Dirichlet) obstacle $A$. Here we follow the discussion in \cite{Colton:2013:IAE}. We assume for simplicity that the obstacle $A$ is a  bounded domain with $C^2$ boundary $\partial A$ (for similar results in the more general case of Lipschitz boundary see \cite[\S 9]{McLean:2000:SES}). The scattering problem problem can be posed as the following exterior Dirichlet problem 
\[
\begin{split}
    \Delta u_s +k^2u_s&=0, \; x \in \mathbb{R}^2 \setminus A\\
    u_s&= -u_i, \; x \in \partial A,
\end{split}
\]
where $u_s$ also satisfies the Sommerfeld radiation condition
\begin{equation*}
    {\displaystyle \lim _{|x|\to \infty }|x|^{\frac {1}{2}}\left({\frac {\partial }{\partial |x|}}-\rmi k\right)u_s(x)=0}
    \label{eq:sommerfeld}
\end{equation*}
where $\partial / \partial |x|$ denotes the radial derivative and the limit is uniform for all directions $x/|x|$ (see \cite[\S 3.4]{Colton:2013:IAE}). The exterior Dirichlet problem has a unique solution $u_s \in H^1_{loc}(\real^2 \setminus A)$ for $\Im(k) \geq 0$ and $u_i|_{\partial \Omega} \in H^{1/2}(\partial A)$, see e.g. \cite[\S 3.2]{Colton:2013:IAE}.  This is clearly the case under the assumptions in \cref{rem:validity}, since $u_i \in C^\infty(\partial\Omega)$.

We seek the scattered field in the form of a mixed single and double-layer potential $\psi \in H^{1/2}(\partial A)$ satisfying
\begin{equation}\label{eqn:scat_field}
u_s(x;k) =\int_{\partial A}dS(y) \bigg(\frac{\partial G}{\partial \nu(y)}(x-y;k)-\rmi\eta G(x-y;k) \bigg) \psi(y)
\end{equation}
where $\eta \neq 0$ with $\eta\Re(k)\geq 0$ is a coupling parameter. This choice guarantees invertibility (see e.g. \cite[\S 3.6]{Colton:1983:IEM}). Here it is necessary to define the corresponding boundary layer operators  for $x\in \partial A$ by
\[\begin{split}
    (S\varphi)(x)&:= 2\int_{\partial A}dS(y) \Mb{G(x-y;k)\varphi(y)},~\text{and}\\
    (K\varphi)(x)&:= 2\int_{\partial A}dS(y) \Mb{\frac{\partial G}{\partial \nu(y)}(x-y;k)\varphi(y)},\\
\end{split} \]
for the single and double layer potential respectively. We note that the operators can be taken as bounded operators  $S,K: L^2(\partial A) \to L^2(\partial A)$ (see e.g. \cite[Theorem 4.4.1]{Nedelec:2001:AEE} for smooth $\partial A$ or \cite[chapter 6]{McLean:2000:SES} for $C^2$  or even Lipschitz $\partial A$).
We also need the following jump relations, letting $z \in \real^2 \setminus \partial A$ 
\begin{equation}\label{eqn:jump}
    \begin{split}
        \lim_{z \to x}\int_{\partial A}dS(y) \Mb{G(z-y;k)\varphi(y)}&= \frac{[S\varphi](x)}{2},~\text{and}\\
        \lim_{z \to x^+}\int_{\partial A}dS(y) \Mb{\frac{\partial G}{\partial \nu(y)}(z-y;k)\varphi(y)}&=\frac{1}{2}\bigg(\varphi(x)+ [K\varphi](x)\bigg),
    \end{split}
\end{equation}
where $z \to x^+$ denotes the limit from the exterior of $A$. Taking the limit of \eqref{eqn:scat_field} as we approach the boundary of $A$ from the exterior and applying \eqref{eqn:jump} yields
\begin{equation}\label{eqn:bint}
    \psi+K\psi-\rmi \eta S\psi = -2u_i|_{\partial A},
\end{equation}
which has a unique solution $\psi$ (see e.g. \cite[\S 3.2]{Colton:2013:IAE}).
We assume that $\partial A$ admits a $2\pi$-periodic parametrization of the form
\[q(\tau) =  (x_1(\tau),x_2(\tau)), \; 0 \leq \tau \leq 2\pi, \]
that is $q([0,2\pi]) = \partial A$ and $q$ is assumed smooth for our numerical experiments. Following \cite[\S 3.5]{Colton:2013:IAE}, one can reformulate \eqref{eqn:bint} as an integral equation of the second kind with a weakly singular kernel. There are several methods to discretize such integral equations, see e.g. \cite{Barnett:2014:HOA} for a review. Here we chose the Kapur-Rokhlin method \cite{Kapur:1997:HOC}, which is based on the trapezoidal rule for periodic functions. In this method, the unknowns are the values of $\psi$ at uniformly spaced points of $[0,2\pi]$. To account for the singularity, the entries in a band of the system matrix are weighted so that the quadrature is exact for polynomials of a given order (6th order in our case).

 We do note that the Kress quadrature \cite{Colton:2013:IAE} was used for in \cite{Guevara:2011:ECA} for computing the scattered fields with $k>0$ and is spectrally accurate. Unfortunately, accuracy of the Kress quadrature degrades for complex $k$. Indeed, the Kress quadrature is obtained by splitting the singular kernel into a singular and non-singular part. The latter requires the evaluation of $J_0(kr)$, which grows exponentially in $\Im (k)$ for fixed $r>0$, see e.g. \cite[\S 10.7]{NIST:DLMF}. The correction weights for the Kapur-Rohklin method only depend on the type of singularity and order of the method. Thus the Kapur-Rokhlin is better adapted for complex $k$.  Convergence for $k$ complex follows from convergence of the method for the real and imaginary parts, considered individually.

\section{Active exterior cloaking for the heat equation}
\label{sec:timedomain}
We now apply the single wavenumber exterior cloaking approach to the time domain heat equation. We recall in \cref{sec:background} other cloaking approaches. We then use the Fourier-Laplace transform to obtain the Helmholtz equation from reasonable heat equation solutions \cref{sec:t2f}. The exterior cloaking approach is applied for different wavenumbers and then put together again in \cref{sec:numexp} via the inverse Laplace transform. The details of the discrete Fourier transform based algorithm we used for this purpose are in \cref{sec:ilt}.
\subsection{Other cloaking approaches for the heat equation}
\label{sec:background}
Cloaking for the heat equation was originally introduced through a change of coordinate system \cite{guenneau2012transformation}, inspired by transformation optics \cite{greenleaf2003anisotropic,pendry2006controlling}. However, this approach leads to an extreme anisotropic thermal conductivity, and even a thermal cloak designed through a regularized geometric transform suffers from limited efficiency in the transient regime \cite{guenneau2012transformation}. A good thermal cloak efficiency requires as many as 10,000 isotropic concentric layers to finely approximate its spatially varying anisotropic conductivity \cite{petiteau2014spectral}. Thus, fabricated metamaterial cloaks with a limited number of layers suffer from reduced efficiency in the transient regime \cite{schittny2013experiments,narayana2012heat,xu2014ultrathin,han2014experimental,hu2018illusion}. For other passive cloaking and mimicking approaches see e.g. \cite{alwakil2017inverse,diatta2010non}. Recent advances in thermal cloaking are thus underpinned by inverse homogenization problems that require heavy computational resources. On the other hand, thermoelectric devices have been proposed to pump the heat flow accurately from one side of a thermal cloak to the other side by adjusting the input current, so that the background temperature field can be restored in a stationary regime \cite{nguyen2015active,huang2020thermal}. In our former work \cite{Cassier:2020:ATC},  we envisioned using Peltier devices (surrounding the object to cloak) to control transient thermal fields generated by a source. Our approach can be viewed as a generalization of that in \cite{xu2019dipole} that considered a single dipole source placed inside the object to cloak. There should be a trade-off between using a single dipole source and numerous monopole and dipole sources to achieve efficient thermal cloaking in the transient regime, which is what motivated the present work. Our analysis is performed in the frequency regime, where we can extend results of \cite{Guevara:2011:ECA} to the Helmholtz equation with complex wavenumbers. Results are then translated in the time domain through inverse Fourier-Laplace transform.

\begin{Rem}
Since our approach is based on the Laplace transform of the time domain heat equation, it is more convenient to assume a zero initial condition. Indeed a non-zero initial condition would appear as a source term for the Helmholtz equation, which would prevent us from using the interior reproduction formula \eqref{eqn:brepfreq}.  However, as noted in \cite{Cassier:2020:ATC}, if the initial condition is a steady state solution to the heat equation (i.e. harmonic) we can use the linearity of the heat equation to apply our approach to $\mathfrak{u}(x,t) - \mathfrak{u}(x,0)$.
\end{Rem}

\begin{Rem}
As we see next, we obtain solutions to the heat equation that achieve exterior cloaking but they can be large as we get close to the new source locations $x_j$. However, as noted in \cite{Cassier:2020:ATC}, it is conceivable to use Peltier devices to physically implement  the interior/exterior reproduction formula for the time domain heat equation \cite{Costabel:1990:BIO}. This procedure allows to replace point-like sources by active surfaces that we call ``extended cloaking devices'', which would keep the temperatures at levels that would be practical to implement.
\end{Rem}

\subsection{From time domain to frequency domain}
\label{sec:t2f}
We now apply our frequency domain cloaking approach to the heat equation. The temperature $\mathfrak{u}(x,t)$ (measured in Kelvin) in a homogeneous isotropic medium satisfies the heat equation,
\begin{equation}\label{eqn:heat}
\rho c \frac{\partial \mathfrak{u}}{\partial t} = \kappa  \Delta \mathfrak{u} +\mathfrak{h},~\text{for }~t > 0,
\end{equation}
where $t$ is the time (s), $\rho$ is the mass density (kg.m$^{-2}$), $c$ is the specific heat (J.K$^{-1}$.kg$^{-1}$) and $\kappa$ is the thermal conductivity (W.K$^{-1}$). Here we assume that $\rho$, $c$ and $\kappa$ are positive constants. The source term is $\mathfrak{h}(x,t)$ (W.m$^{-2}$) and assumed causal, i.e. $\mathfrak{h}(x,t) = 0$ for $t<0$. 
For simplicity we assume a zero initial condition and consider 
\begin{equation}\label{eqn:heat_bis}
\frac{\partial \mathfrak{u}}{\partial t} = \sigma  \Delta \mathfrak{u} +\frac{\mathfrak{h}}{\rho c},~\text{for }~t > 0,
\end{equation}
 where $\sigma=\kappa/\rho c$ is the thermal diffusivity ($m^2$.s$^{-1}$).
Assuming further that the source term $\mathfrak{h}(x,t)$ satisfies the growth condition \eqref{eq:growth} with $L^2(\real^2)$ norm, $\alpha \geq 0$ and $p=0$, we can see that $e^{-\xi t}\mathfrak{h}(x,t) \in L^2([0,\infty),L^2(\real^2))$ for any $\xi > \alpha $. Using \cite[Corollary 2, p238]{Dautray:1992:MAN} it is possible to conclude that \eqref{eqn:heat_bis} admits a unique solution $\mathfrak{u}(x,t)$ satisfying $e^{-\xi t}\mathfrak{u}(x,t) \in L^2([0,\infty),H^2(\real^2))$ for any $\xi > \alpha$. This allows to define the Fourier-Laplace transform \eqref{eqn:lp_transform} of all terms in \eqref{eqn:heat_bis} on the half plane $\complex_\alpha^+$, thus obtaining  the Helmholtz equation
\begin{equation} \label{eqn:im_helm}
\Delta u(x;\omega)+ \frac{\rmi\omega}{\sigma}u(x; \omega) = -\frac{h(x,\omega)}{\rho c},
\end{equation}
where the wavenumber is $k = \rmi \sqrt{-\rmi\omega/\sigma}$, using the principal value of the square root and $h$ is the Fourier-Laplace transform of  $\mathfrak{h}$. We note that $\Re(k^2) = \Re(\rmi\omega/\sigma) < 0$, whenever $\Im(\omega>0)$,  which guarantees that the Helmholtz equation satisfies a form of the maximum principle for any $\omega \in \complex^+_\alpha$ (since $\alpha \geq 0$), for $x$ outside of the support of the source $h$ (see \cref{sec:stability} and \eqref{eq:calpha} for the definition of $\complex_\alpha^+$).

\subsection{Numerical experiments}
\label{sec:numexp}
\begin{figure}
    \centering
    \begin{tabular}{cc}
    \includegraphics[width=0.45\textwidth]{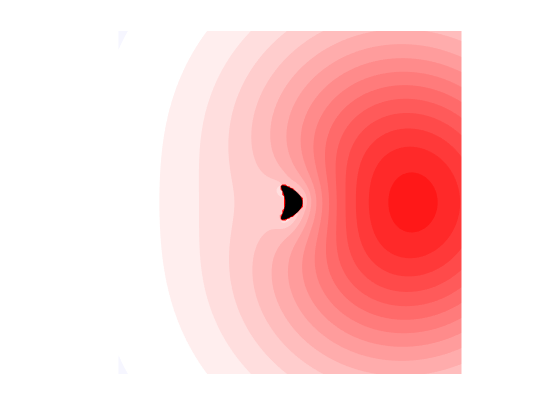} & \includegraphics[width=0.45\textwidth]{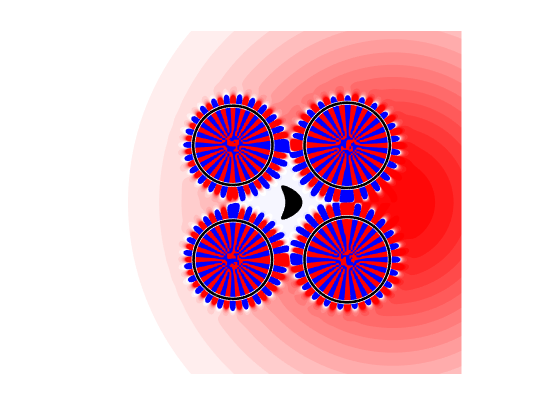}\\
    (a) & (b)
    \end{tabular}
    \caption{(a) Temperature distribution at time $t=4$ with $\sigma=1.5$ (see \eqref{eqn:heat_bis}) resulting from a point source at location $(8,5)$ in the presence of a ``kite'' object with homogeneous Dirichlet boundary condition. (b) The same object and source as in (a), but with the cloaking devices activated. The temperatures outside of the black circles are bounded by $\mathfrak{u}_{max} \approx 6.2$ (see \eqref{eq:umax}) for $t \in [0,4]$. The computational domain was $[0,10]^2$. The linear color scale spans temperatures in $[-0.0133,0.0133]$ and range from blue (negative) to red (positive), with zero been represented in white. [See also movie in supplementary material]}
    \label{fig:time_dom}
\end{figure}
We show in \cref{fig:time_dom} a numerical simulation of active exterior cloaking of a Dirichlet object (a ``kite'' with constant zero temperature at its boundary) and compare it to the case where there is no cloaking devices. We used $\sigma = 1.5$ in \eqref{eqn:heat_bis}. As can be seen from the time snapshot in \cref{fig:time_dom}, the isotherm lines without the cloaking devices are significantly different from those of the point source that we used as the incident field $u_i$, this is because of the field ``scattered'' by the object. For an observer far from the cloaking devices, the isotherms appear consistent with those of a point source, so it is hard for the observer to detect the object from thermal measurements. In our numerical experiments, the region $\Omega$ is a disk centered at $(5, 5)$ and with radius $\delta_C = 10/6$ enclosing the ``kite'' object. We moved the distribution of monopoles and dipoles to 4 new source locations determined as in \cref{fig:cart_config} with $\delta_D = 5\sqrt{2}/3$. We note that three new source locations would have been sufficient, as in \cite{Guevara:2011:ECA}. The fields are calculated on $[0,10]^2$ using a $200 \times 200$ uniform grid. The incident field is generated by a point source at $y=(8,5)$. The integral in \cref{thm:conv} is approximated with 256 uniformly placed points and the series in \eqref{eqn:ext_trunc} uses the truncation $M=22$. The boundary of the scatterer is discretized using 512 equally spaced points on the parametric representation of $\partial A$ and the scattered fields (in the frequency domain) are calculated according to the scheme in \cref{sec:numericsfreq}. The frequency domain calculation is performed for 2050 frequencies and the Laplace transform is inverted using a Fast Fourier Transform based method (see \cref{sec:ilt}). 

Because the multipolar sources are singular at the $x_j$, the cloaking field diverges as we approach the $x_j$. This could limit the physical implementation (e.g. because the material starts degrading with such high temperatures). Of course, we may use the Green exterior representation formula (see e.g. \cite{Costabel:1990:BIO,Cassier:2020:ATC}) to replace each of the multipolar sources by a monopole and dipole distribution on surface containing the multipolar source. These active surfaces  or ``extended cloaking device'' could be realized in practice using Peltier devices \cite{nguyen2015active} and the temperatures do not need be unreasonably large. To illustrate this, we display in \cref{fig:time_dom} black circles centered at the new source positions, outside of which we are guaranteed to have $|\mathfrak{u}_e^{(M)}(x,t)| \leq   \mathfrak{u}_{max}$. Our choice $\mathfrak{u}_{max}$ is  
\begin{equation}
 \label{eq:umax}
 \mathfrak{u}_{max} = 100\max_{(x,t) \in \Omega \times [0,T] } |\mathfrak{u}_i(x,t)|. 
\end{equation}
This choice is not a statement of what is feasible, but simply for illustration purposes. In \cref{fig:time_dom}, we have $T=2$ and $\mathfrak{u}_{max} \approx 6.2$. 

\begin{figure}
	\centering
	\includegraphics[width=0.5\textwidth]{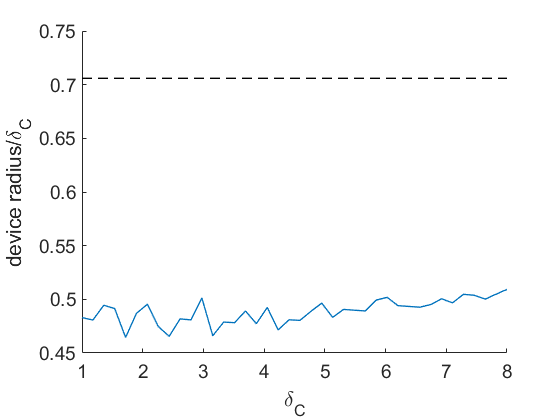}
	\caption{The blue curve corresponds to circular active surface radii (relative to the scaling parameter $\delta_C$) to achieve exterior cloaking of a point source located at $(10,1)$ for $40$ different values of the scaling parameter $\delta_C$ and for $\Omega$ centered at $(10,10)$ (see \cref{fig:cart_config}). The dotted line corresponds to the value for which the circles touch. Since the blue data points are below the dotted line, there are gaps between the circular active surfaces, showing that even the ``extended cloaking devices'' do not completely surround the object (exterior cloaking).
	}
	\label{fig:circles}
\end{figure}

To show that we are achieving exterior cloaking even when replacing the multipolar sources by extended cloaking devices (circular active surfaces), we changed the scale $\delta_C$ (radius of $\Omega$) of the cloaking configuration in \cref{fig:cart_config} with $\delta_D =\sqrt{2} \delta_C$, keeping $\Omega$ as a disk with a fixed center $(10,10)$ and a fixed point source positioned at $(10,1)$, which generates the incident field (see also \cref{sec:evaluation}).
For all values of $\delta_C$, we used the same diffusivity $\sigma = 1.3$ in \eqref{eqn:heat_bis}, truncation $M=22$ and 128 points to discretize a parametric representation of $\partial\Omega$. The temperature fields where evaluated  using 130 wavenumbers. The computation was repeated for $40$ equally spaced $\delta_C$ in the interval $[1,8]$, chosen so that the divergence region $R$ from \cref{thm:conv} does not include the source location. For each $\delta_C$, the cloak field was evaluated on a $100 \times 100$ uniform grid of the square $[0,20]^2$ and the circles outside of which $|\mathfrak{u}(x,t)| \leq \mathfrak{u}_{max}$ were determined with $T = 1$ in \eqref{eq:umax}. We display in \cref{fig:circles} the radius of these circles relative to $\delta_C$ and as a function of $\delta_C$. The dotted line in \cref{fig:circles} corresponds to the radius for which the circular active surfaces would touch and match the divergence region $R$. As we can see from \cref{fig:circles}, the circular active sources do not touch, and thus we have exterior cloaking even in this situation. Roughly speaking, according to \cref{fig:circles}, the ``urchins'' have a radius that is about 70\% of the radius of the gray circles in \cref{fig:cart_config}.

\subsection{From the frequency domain to the time domain}\label{sec:ilt}
For convenience we express the numerical algorithm we use to go from frequency domain to time domain in terms of the Laplace transform of $\mathfrak{u}(x,t)$ rather than the Fourier-Laplace transform \eqref{eqn:lp_transform}. Under the same growth condition assumption \eqref{eq:growth} on $\mathfrak{u}(x,t)$, its Laplace transform is
\begin{equation}
    \sfu(x,s) = \int_0^\infty dt [e^{-st} \mathfrak{u}(x,t)],
 \label{eq:laplace}
\end{equation}
which is well defined for $\Re(s) > \alpha$, where  $\alpha \geq 0$ is defined in section \ref{sec:t2f}. Clearly we have $\sfu(x,s) = u(x,\rmi s)$, where the right hand side is the Fourier-Laplace transform of $\mathfrak{u}$ given in \eqref{eqn:lp_transform}. The inverse Laplace transform is then given by
\begin{equation}\label{eqn:inv_lp}
    \mathfrak{u}(x,t)= \frac{1}{2\pi \rmi} \int_{c-\rmi\infty}^{c+\rmi\infty}ds [e^{s t}\sfu(x,s)],
\end{equation}
for any $c>\alpha$.
\footnote{We note that the assumptions of \cref{sec:t2f} (i. e. the growth control of the source term,  the zero initial condition and the causality of $\mathfrak{u}$) imply that for any $\xi > \alpha \geq 0$ : $e^{-\xi t}\mathfrak{u}(x,t) \in L^2([0,\infty),H^2(\real^2))$ and (using the heat equation) that  $e^{-\xi t} \partial_t \mathfrak{u}\in L^2([0,\infty),L^2(\real^2))$. 
Thus, one has that $\sfu(\cdot,s)= s^{-1}\, \mathcal{L}(\partial_t \mathfrak{u})(\cdot,s)$ is analytic with respect to $s$ for $\Re (s) > \alpha$, where $\mathcal{L}$ stands for the Laplace transform. Since $\|\mathcal{L}(\partial_t \mathfrak{u})(\cdot,s)\| \to 0$ for $\Re(s) \geq c> \alpha$ and $|s| \to \infty$, we get $\|u(\cdot,s)\|=o(|s|^{-1})$. We need a little more of decay to use the formula \eqref{eqn:inv_lp}  for the inverse Laplace transform. It is enough to assume that there exists $\varepsilon>0$ such that $\|u(\cdot,s)\| \leq C |s|^{-(1+\varepsilon)}$ for $\Re(s) \geq c> \alpha$. Then,  \eqref{eqn:inv_lp} is well defined as a Bochner integral valued in $L^2(\mathbb{R}^2)$ (see e.g. the proof of \cite[Theorem 2.5.1]{arendt2001cauchy} or \cite{Sayas:2016:RPT}) and thus also pointwise for a. e. $x\in\mathbb{R}^2$ and all $t \geq 0$.    
}
We follow the numerical method in \cite{HSU:1987:NIC} for computing the inverse Laplace transform by approximating it with a discrete Fourier transform (or DFT, which can be evaluated efficiently with the Fast Fourier Transform or FFT see e.g. \cite{Frigo:2005:DIF} for definition). For a review of numerical inverse Laplace transform methods see e.g. \cite{Cohen:2007:NML}. The idea is that to a uniform grid of the time interval $[0,T]$ with $N$ points, i.e. $t_p=p\Delta t \in [0,T]$, $p=0,\ldots,N-1$, $\Delta t =T/N$, we associate the discretization of a dual variable $w$  given by $w_q = q \Delta w$, $q=0,\ldots,N-1$, $\Delta w = 2\pi/T$. The discretizations are chosen such that $t_p w_q = 2\pi p q / N$, which is the negative of the phase of the complex exponential in the DFT of length $N$. Using the change of variables $s = c+\rmi w$ in \cref{eqn:inv_lp} and approximating with a Riemann sum on a finite interval yields the following (we assume $\mathfrak{u}(x,t)$ is real)
\[ 
\begin{split}
\mathfrak{u}(x,t) &= \frac{1}{2 \pi} \int_{-\infty}^\infty dw \Mb{\exp(ct)\exp(\rmi w t)\sfu(x,c+\rmi w)}\\
&\approx \frac{1}{2 \pi} \sum_{q = -(N-1)}^{N-1}\Delta w \Mb{\exp(ct)\exp(\rmi w_qt)\sfu(x,c+\rmi w_q)}\\
&= \frac{\exp(ct)}{T} \sum_{q = -(N-1) }^{N-1}  \Mb{\exp\M{\rmi w_qt }\sfu\M{x,c+\rmi w_q}}\\
&=2\frac{\exp(ct)}{ T}  \Re \Mb{ \sum_{q = 0 }^{N-1}  \Mb{\exp\M{-\rmi w_qt }\sfu\M{x,c-\rmi w_q}}-\frac{\sfu(x,c)}{2}},
\end{split}
\]
where the last equality follows by the the symmetry for Laplace transforms of real functions ($\overline{\sfu(x,s)} = \sfu(x,\overline{s})$). By evaluating at the $t_p$, our approximation can be written in terms of the real part of a DFT:
\[ 
\begin{split}
\mathfrak{u}(x,t_p)&=2\frac{\exp(ct_p)}{T} \Re\Mb{\sum_{q = 0}^{N-1}  \Mb{\exp\M{-\rmi t_pw_q}\sfu\M{x,c-\rmi w_q}}-\frac{\sfu(x,c)}{2}}\\
&=2\frac{\exp(ct_p)}{T} \Re\Mb{\texttt{fft}(\sfu\M{x,c-\rmi w_q},q=0,\ldots,N-1)-\frac{\sfu(x,c)}{2}}.
\end{split}
\]
Here $\texttt{fft}(v)$ represents the DFT of a vector $v$ of length $N$, as defined in \cite{Frigo:2005:DIF}. Thus we end up evaluating the frequencies $\omega_q = w_q + \rmi c$, $q=0,\ldots,N-1$.
As noted in \cite{HSU:1987:NIC}, the convergence only holds for $[0, T/2)$, so in practice we use $\widetilde{N}=2N+2$ and only use the first $N+1$ time steps to give convergence on the interval $[0,T]$ . Here the additional two frequencies mean the terminal time step is $T$ as opposed to $T -\Delta T.$ Since the heat equation solutions we consider decay, we take $\alpha=0$ and set $c = \alpha- \frac{\Delta w}{2 \pi }\ln(10^{-6}) > 0$. We note that there are methods to speed up convergence of this class of numerical inverse Laplace transform, see e.g. \cite{Brancick:2002:MOM}.

\section{Summary and perspectives}
\label{sec:discussion}
In our earlier work on active exterior cloaking for the parabolic heat equation \cite{Cassier:2020:ATC}, we noted in the concluding remarks that our approach could also be tied to the active exterior cloaking strategies for the Helmholtz equation by going
to Fourier or Laplace domain in time. Here we have showed that it is possible to cloak objects from thermal measurements by using active heat sources, starting from a zero temperature condition. We believe that our work opens up a new path for active cloaking for a variety of physical situations, in addition to the class of differential equations of the form \eqref{eq:poly}.  Indeed,
complex wavenumbers make it possible to model pseudodifferential operators in time such as fractional time derivatives \cite{podlubnv1999fractional} and integro-differential equations. This opens new avenues in active exterior cloaking. To give an example of this flexibility, consider a (non-dimensionalized) heat equation with memory that arises when considering homogenized diffusion models in fractured media \cite{Hornung:1990:DMF}
\begin{equation}
    \frac{\partial \mathfrak{u}}{\partial t} + \int_0^t \mathfrak{p}(t-\tau) \frac{ \partial\mathfrak{u}}{\partial t}(\tau) d\tau = \Delta \mathfrak{u} +  \mathfrak{h},
\end{equation}
where $\mathfrak{h}(x,t)$ is a source term and $\mathfrak{p}(t)$ is a convex monotone decreasing history function with a singularity at $t=0$ (e.g. $\mathfrak{p}(t) = \alpha \exp[- \alpha t]$).
Indeed by taking Fourier-Laplace transform and assuming zero-initial conditions we get in the frequency domain 
\begin{equation}
 \Delta u(x,\omega) -(-\rmi\omega + p(\omega)) u(x,\omega)   = - h(x,\omega).
\end{equation}

The present work allows us to also study active cloaking in the context of diffusive photon density waves governed by
\begin{equation}\label{eqn:photon}
\frac{\partial \mathfrak{u}}{\partial t} + \mu \mathfrak{u}= \sigma  \Delta \mathfrak{u} +\mathfrak{h},~\text{for }~t > 0,
\end{equation}
where $\mu>0$ is an absorption coefficient, $\sigma$ is a conductivity and $\mathfrak{h}$ represents the photon current density (photon flow per unit surface and per unit time).
Making use of the Fourier-Laplace transform \eqref{eqn:lp_transform} and assuming zero initial conditions, \eqref{eqn:photon} takes the form 
of the Helmholtz equation in the context of diffusion wave scattering
\begin{equation}\label{eqn:lp_helmholtz_photon}
\Delta u(x;\omega)+ \frac{\rmi\omega-\mu}{\sigma}u(x; \omega) = -\frac{h(x,\omega)}{\sigma},
\end{equation} 
where we note that \eqref{eqn:lp_helmholtz_photon} reduces to \eqref{eqn:im_helm}, when $\omega\gg \mu$.
Scattering cancellation of such diffusive waves has been addressed in \cite{schittny2014invisibility,farhat2016cloaking}.

Moreover, advection-diffusion problems play a prominent role notably in diffusion and mixing of fluid flow modelled by \cite{constantin2008diffusion}
\begin{equation}\label{eqn:heatdrift}
\rho c \frac{\partial \mathfrak{u}}{\partial t} = - v\cdot\nabla u + \kappa  \Delta \mathfrak{u} +\mathfrak{h},~\text{for }~t > 0,
\end{equation}
where $v$ is a constant velocity, $\rho$ is a mass density, $\kappa$ a conductivity and $c$ the heat capacity.  The same equation is known as the Fokker-Planck equation and is central to models for transport of salt, heat, buoys, and markers in geophysical flows \cite{murphy2020spectral,coti2020homogenization}.
It turns out that one can recast \eqref{eqn:heatdrift} using the exponential variable transform
\cite{LI1991255}
\begin{equation}
\mathfrak{u}(x,t)=\exp[ \frac{v}{2\, \kappa}\cdot x] \mathfrak{w}(x,t) \quad \mbox{ and } \quad \mathfrak{h}(x,t)=\exp[ \frac{v}{2\, \kappa}\cdot x] \, \mathfrak{g}(x,t)
\end{equation}
together with the Fourier-Laplace transform \eqref{eqn:lp_transform}, into the Helmholtz equation (assuming zero initial conditions)
\begin{equation} 
\Delta w(x;\omega)-\tau^2 w(x; \omega) = -\frac{
g(x,\omega)}{\kappa},
\end{equation}
where
$\tau^2 =(\mid{v}\mid/2\,\kappa)^2-\rmi\omega/\sigma$, and $\sigma=\kappa/(\rho c)$.

Finally, the exterior cloaking theory which we developed may allow us to also achieve exterior cloaking in the context of Maxwell-Cattaneo heat waves governed by \cite{RevModPhys.61.41}
\begin{equation}\label{eqn:maxcat}
\tau\frac{\partial^2 \mathfrak{u}}{\partial t^2}+\frac{\partial \mathfrak{u}}{\partial t}= \kappa  \Delta \mathfrak{u} +\tau\sigma\Delta\frac{\partial\mathfrak{u}}{\partial t} +\mathfrak{h},~\text{for }~t > 0,
\end{equation}
where $\kappa$ is the thermal conductivity, $\sigma$ accounts for diffusive phenomena, $\tau$ is the thermal relaxation time (that corresponds
to the time it takes for a medium to reduce its temperature to half). We assume that $\kappa$, $\tau$, $\sigma$ are positive constants.
Making use of the Fourier-Laplace transform \eqref{eqn:lp_transform} with zero initial conditions, \eqref{eqn:maxcat} takes the form 
of the Helmholtz equation in the context of diffusion wave scattering
\begin{equation}\label{eqn:lp_helmholtz_photon_maxwell}
\Delta u(x;\omega)+ \frac{\omega(\omega\tau+\rmi)}{\kappa-\rmi\omega\tau\sigma}u(x; \omega) = -\frac{h(x,\omega)}{\kappa - \rmi \omega \tau \sigma},
\end{equation} 
and thus one can define $k^2 =\omega(\omega\tau+\rmi) / (\kappa-\rmi\omega\tau\sigma)$ that unlike for the Fourier heat equation \eqref{eqn:im_helm} can lead to propagating features (when $\Re(k^2)>0$).
Scattering cancellation of such diffusive waves has been addressed in \cite{farhat2019scattering}.

Another situation where complex wavenumbers may play a prominent role is for
in-plane pressure and shear elastodynamic waves propagating in passive, dissipative, active or even viscoelastic media (the case of anti-plane shear waves would be covered by the 2D scalar Helmholtz equation with complex wavenumber). The latter, viscoelastic media, would require using the Helmholtz decomposition proposed in \cite{bretin2011green} $ u=\nabla\Phi+\nabla\times\Psi$, with $\Psi$ a divergence free vector field, in the vector Navier equation
\begin{equation}
(\tilde{\lambda}+2\tilde{\mu})\nabla\nabla\cdot u+\tilde{\mu}\nabla^2 u+\rho\omega^2 u= 0
\end{equation} 
where $\tilde{\lambda}=\lambda+\eta_p{\mathfrak{M}}$, $\tilde{\mu}=\mu+\eta_s{\mathfrak{M}}$, $\mathfrak{M}$ being a convolution operator with certain power law (named after Szabo and Wu \cite{szabo2000model}), $\lambda$, $\mu$ are the usual Lam\'e parameters, $\rho$ is the density and $\eta_s \; , \; \eta_p \ll 1$. Our approach would then be applied to Helmholtz equations with
the complex wavenumbers for shear (s) and pressure (p) waves
\begin{equation}
k_m^2(\omega)=\omega^2\M{1-\frac{\nu_m}{c_m^2}\mathcal{M}(\omega)} \; , \; m=s, \; p \; ,
\end{equation}
with $c_p=\sqrt{(\lambda+2\mu)/\rho}$, $c_s=\sqrt{\mu/\rho}$, the pressure and shear wave velocities, respectively, $\nu_p=(\eta_p+2\eta_s)/\rho$, $\nu_s=\eta_s/\rho$ and where the multiplication operator $\mathcal{M}$ is the Fourier transform of the convolution operator ${\mathfrak{M}}$. The Hodge decomposition can also be applied to split the elastodynamic wave equation with isotropic viscoelasticity into two acoustic wave equations in time (P and S waves) that can be transformed into the Helmholtz equation, see e.g. \cite{Albella:2018:SLI}.

There may also be ways of adapting our approach to other differential operators in space. One example would be to consider constant anisotropic media (e.g. coming from a homogenization approach). We believe that our strategy could be also applied to \eqref{eq:poly} wherein the Laplacian is replaced by a bi-Laplacian in the right-hand side. This problem arises when modeling flexural waves in thin elastic plates, which are governed by the bi-harmonic (Kirchhoff-Love) equation. Active cloaking in this context has been considered in \cite{Oneil:2015:ACI}. We are also considering applying our theory to active cloaking for flexural gravity waves in floating thin elastic plates that would involve a tri-Laplacian in the right-hand side of \eqref{eq:poly} \cite{farhat2020scattering}.

We conjecture that our approach can be adapted to the 3D Helmholtz equation for complex wavenumbers, which would allow us to address cloaking in general dispersive media (including media with losses or gain). 
Finally we note that open questions related to gain media ($\Im(k)<0$) remain. In particular what is a sensible functional space setting for the exterior Green representation formula and for the scattering problem in gain media.
We believe that the approach of exterior cloaking which we have developed in this article allows us to cover a broad range of problems for active cloaking of diffusion and wave phenomena. 

\appendix
\section{Proof of \cref{lem:gbound}}
\label{app:gbound}

\begin{proof} 
 {\it Step 1:} We prove first  \eqref{eq.Jn}.
The entire function $J_n$  is defined \cite{Watson:1995:TTB}[Chapter, formula 1 p. 15] via the following power series:
\begin{equation}\label{eq.defJn}
J_n(z):= \sum_{k=0}^{\infty} \frac{(-1)^k }{ k! \, (n+k)!}  \Big(\frac{z}{2}\Big)^{n+2k}=  \frac{1}{n!}\Big(\frac{z}{2}\Big)^{n}+\sum_{k=1}^{\infty} \frac{(-1)^k }{ k! \, (n+k)!} \Big(\frac{z}{2}\Big)^{n+2k}, \quad \, \forall z\in \mathbb{C}.
\end{equation}
As $1/(n+k)!\leq 1/(n+1)!$ for $k\geq 1$, one gets
 \begin{equation}\label{eq.estimreminderJn}
 \Big|  \sum_{k=1}^{\infty} \frac{(-1)^k }{ k! \, (n+k)!} \Big(\frac{z}{2}\Big)^{n+2k}\Big| \leq  \frac{1}{(n+1)!}\Big(\frac{|z|}{2}\Big)^{n+2} \Big[ \sum_{k=1}^{\infty} \frac{1}{k!}\Big(\frac{|z|}{2}\Big)^{2k-2} \Big] 
\leq \frac{C_{K_1}} {(n+1)!}\Big( \frac{|z|}{2}\Big)^{n+2}
 \end{equation}
 where the positive constant $C_{K_1}$ is defined by 
 \begin{equation}\label{eq.defCK1}
C_{K_1}=\max_{z\in K_1} f(|z|) \ \mbox{with $f$ entire defined by} \ f(z) :=  \sum_{k=1}^{\infty} \frac{1}{k!}\big(\frac{z}{2}\big)^{2k-2},\ \ \forall z\in \mathbb{C}. 
\end{equation}
We point out that $f(z)=4(\exp(z^2/4)-1)/ z^2$ for $z\neq 0$ and $f(0)=1$. Thus,
\eqref{eq.Jn} is an immediate consequence of \eqref{eq.defJn} and  \eqref{eq.estimreminderJn}. \\[12pt]
\noindent {\it Step 2:} The asymptotic formula \eqref{eq.Jnprime} for  $J_n'$ follows  immediately from the asymptotic expansion \eqref{eq.Jn} and the recurrence formula  (see \cite{Watson:1995:TTB} formula 2 page 45)
 $
 J_{n}'(z)= \frac{1}{2}( J_{n-1}(z))- J_{n+1}(z)).
 $
 An other way  to obtain the  formula  \eqref{eq.Jnprime}  is to derive the power series that defines $J_n$ and applies the same method as for formula \eqref{eq.Jn}. 
 \\[12pt]
{\it Step 3:} We now  prove \eqref{eq.Hn}.  By definition of the Hankel function $H_n^{(1)}$ (see \cite{Watson:1995:TTB}, formula (1) page 73) , one has for all $z\in \mathbb{C}\setminus (-\infty,0]$:
\begin{equation*}
H_n^{(1)}(z)- \frac{\rmi  \, (n-1)! }{\pi} \Big(\frac{2}{z}\Big)^n =J_n(z)+\rmi Y_n(z)-\frac{\rmi  \, (n-1)! }{\pi} \Big(\frac{2}{z}\Big)^n  
\end{equation*}
Using \eqref{eq.defJn} and the series representation of $Y_n$ on $\mathbb{C}\setminus (-\infty,0]$ (see  \cite{Watson:1995:TTB}, formula (3) page 62 and (2) page 64 or formula (10.8.1) of \cite{NIST:DLMF})  on the previous expression leads to:
\begin{equation}\label{eq.Hnbis}
 H_n^{(1)}(z)- \frac{\rmi  \, (n-1)! }{\pi} \Big(\frac{2}{z}\Big)^n  = \frac{\rmi}{\pi} \Big(\frac{2}{z}\Big)^{n} \sum_{k=1}^{n-1} \frac{(n-k-1)!}{k!} \big(\frac{z}{2}\big)^{2k}+ I_2(n,z) 
\end{equation}
where
$$
I_2(n,z)=\big(\frac {z}{2}\big)^{n}\Big(\sum_{k=0}^{\infty} \Big(\frac{\rmi}{\pi}\big(2 \, \ln(z/2)-\psi(k+1)-\psi(k+n+1)\big)+1\Big) \frac{(-1)^k}{k! (n+k)!} \Big(\frac{z}{2}\Big)^{2k},
$$
where $\ln$ is the principal value of the logarithm function with a branch cut on  $(-\infty,0]$ and  $\psi:=\Gamma'/\Gamma$ with $\Gamma$ the well-known Gamma function.  
We estimate now the two terms appearing in  \eqref{eq.Hnbis}. For the first one, one obtains 
\begin{eqnarray}\label{eq.HnTer}
\Big|\frac{\rmi}{\pi} \Big(\frac{2}{z}\Big)^{n} \sum_{k=1}^{n-1} \frac{(n-k-1)!}{k!} \big(\frac{z}{2}\big)^{2k}\Big| &\leq &\frac{(n-2)!}{\pi}   \Big(\frac{2 }{|z|}\Big)^{n-2} \Big[ \sum_{k=1}^{n-1} \frac{1}{k!}\Big(\frac{|z|}{2}\Big)^{2k-2} \Big]. \nonumber \\
&\leq  &\frac{C_{K_2}}{\pi} (n-2)!  \Big(\frac{2 }{|z|}\Big)^{n-2},
\end{eqnarray}
with the constant $C_{K_2}$ defined by replacing the compact  $K_1$ by $K_2$ in  \eqref{eq.defCK1}.
As the function $\psi$ evaluated on integers is given by  (see \cite[\S 5.4.14]{NIST:DLMF}  formula 5.4.14)
$$
\psi(m+1):= \sum_{p=1}^{m} \frac{1}{p} +\gamma \mbox{ for } m\geq 1 \ \mbox{ with }\ \gamma=:\Gamma'(1)=\psi(1) \mbox{ the Euler constant,}
$$ 
the second term of \eqref{eq.Hnbis} can be bounded by
\begin{multline*}
|I_2(n,z)|\leq
  \frac{1}{\pi} \frac{(n-2)!}{(n-2)!}  \Big(\frac{2 }{|z|}\Big)^{n-4}   \Big(\frac {|z|}{2}\Big)^{2n-4} \times\\ \Big[\sum_{k=0}^{\infty}\Big( 2  |\, \ln\big(\frac{z}{2}\big)|+C+  \sum_{p=1}^{k+1} \frac{1}{p}+  \sum_{p=1}^{n+k} \frac{1}{p}\Big)\frac{1}{k! (n+k)!} \big(\frac{|z|}{2}\big)^{2k}   \Big],
\end{multline*}
with $C=2\gamma+1$ (we point out that we use the inequality $\psi(k+1)\leq\psi(k+2)$ for $k\geq 0$ to avoid the particularity of the case $k=0$). Then, using the following inequality $\sum_{p=1}^{m} 1/p\leq \ln(m)+1$ for $m\geq 1$ (obtained by comparison of the harmonic series with the integral) and the fact that $\ln(k+1)\leq \ln(n+k)$, one gets that 
\begin{multline*}
 |I_2(n,z)|\leq  \frac{1}{\pi}  \frac{(n-2)!}{(n-2)!}  \,  \Big(\frac{2 }{|z|}\Big)^{n-2}  \Big(\frac {|z|}{2}\Big)^{2n-2}
 \times\\\Big[\sum_{k=0}^{\infty}\big( 2 |\, \ln\big(\frac{z}{2}\big)|+C+2\ln(n+k)+2 \big)\frac{1}{k! (n+k)!} \big(\frac{|z|}{2}\big)^{2k}   \Big] .
\end{multline*}
Notice that 
$$
 \frac{\ln(n+k)}{(n+k)!} \leq \frac{\ln(n+k)}{ (n+k) (n+k-1)!} \leq \frac{1}{(n+k-1)!} \leq \frac{1}{(k+1)!}  \ \mbox{and } \ \frac{1}{(n+k)!}\leq  \frac{1}{(k+1)!}     \mbox{ for $n\geq 2$},
 $$
 so we obtain that
 $$
 |I_2(n,z)|\leq  \frac{1}{\pi}  \frac{(n-2)!}{(n-2)!}  \,  \Big(\frac{2 }{|z|}\Big)^{n-2}  \Big(\frac {|z|}{2}\Big)^{2n-2}  \big( 2 |\, \ln\big(\frac{z}{2}\big)|+C+4 \big)  \Big[\sum_{k=0}^{\infty} \frac{1}{k! (k+1)!}  \big(\frac{|z|}{2}\big)^{2k}   \Big] .
 $$
 By introducing  the entire functions,
 $$
g(z)=\sum_{n=2}^{\infty}   \frac{1}{(n-2)!} \Big(\frac {|z|}{2}\Big)^{2n-2} \ \mbox{ and } \ h(z)=\sum_{k=0}^{\infty} \frac{1}{k! (k+1)!}  \big(\frac{z}{2}\big)^{2k}, \quad \forall z\in \mathbb{C},
 $$
it follows that 
 $$
 |I_2(n,z)|\leq \tilde{C}_{K_2} (n-2)!  \Big(\frac{2 }{|z|}\Big)^{n-2} \mbox{ with }  \tilde{C}_{K_2}=\frac{1}{\pi}  \max_{K_2}\,g(|z|) (2 |\ln\big(\frac{z}{2}\big)| +C+4) h(|z|).
 $$
\end{proof}

{\bf Data Access.} We provide the following supplementary materials. (i) A movie animating \cref{fig:time_dom}. (ii) The Matlab code to reproduce the figures \cref{fig:freq_err,fig:space_err,fig:rep,fig:freqvtrunc,fig:max,fig:time_dom,fig:circles} is available in the repository \cite{thecode}.

{\bf Author Contributions.} SG initiated the project. MC proved \cref{lem:gbound}. MC and TD proved the error estimates for the Graf addition formula. MC, TD and FGV proved the error estimates for the truncated cloaking field. TD and FGV performed the numerical experiments. SG found that the approach was valid even for media with gain. All authors contributed to the writing. All authors contributed other idea applications. All authors gave final approval for publication and agree to be held accountable for the work performed therein.

{\bf Funding.} TD and FGV were supported by the National Science Foundation Grant DMS-2008610.
\bibliographystyle{abbrv}
\bibliography{extcloak}
\end{document}